\documentclass[journal]{IEEEtran}
\usepackage{nicefrac}

\usepackage{array}

\usepackage{times}
\usepackage{enumerate}

\usepackage{graphicx}
\usepackage[usenames]{color}

\usepackage{amsmath}  
\usepackage{amssymb}  
\usepackage{mathrsfs} 

\usepackage{theorem}  
\usepackage{cite}     

\usepackage{comment}  

\usepackage{amscd}
\usepackage{latexsym}

\usepackage{upref}
\usepackage{amsfonts}

\definecolor{lgray}{gray}{0.4}

\usepackage{pict2e}
\usepackage[caption=false,font=footnotesize]{subfig}

\newcommand{\cF}{{\cal F}}

\newcommand{\bR}{\mathbb{R}}
\newcommand{\bZ}{\mathbb{Z}}

\newcommand{\DD}{\mathrm{DD}}
\newcommand{\oDD}{\overline{\mathrm{DD}}}

\theoremstyle{plain}
\theorembodyfont{\normalfont\slshape}

\newtheorem{thm}{Theorem$\!$}
\newenvironment{theorem}{\begin{thm}\hspace*{-1ex}{\bf.}}{\end{thm}}

\newtheorem{prop}[thm]{Proposition$\!$}

\newtheorem{lem}[thm]{Lemma$\!$}
\newenvironment{lemma}{\begin{lem}\hspace*{-1ex}{\bf.}}{\end{lem}}

\newtheorem{cor}[thm]{Corollary$\!$}
\newenvironment{corollary}{\begin{cor}\hspace*{-1ex}{\bf.}}{\end{cor}}

\newtheorem{scheme}{Scheme$\!$}

\newtheorem{construction}{Construction$\!$}

\newtheorem{defi}{Definition$\!$}
\newenvironment{definition}{\begin{defi}\hspace*{-1ex}{\bf.}}{\end{defi}}

\newtheorem{exam}{Example$\!$}
\newenvironment{example}{\begin{exam}\hspace*{-1ex}{\bf .}}{\end{exam}}

\begin{document}
\title{Distinct Difference Configurations: Multihop Paths and Key Predistribution in Sensor Networks}
\author{Simon~R.~Blackburn, Tuvi~Etzion, Keith~M.~Martin and Maura~B.~Paterson%
\thanks{This work was supported in part by EPSRC grants EP/D053285/1 and EP/F056486/1, and Israel Science Foundation grant 230/08.}%
\thanks{S.R.~Blackburn, K.M.~Martin and M.B.~Paterson are with the Department of Mathematics, Royal Holloway, University of London, Egham, Surrey TW20 0EX.  T.~Etzion is with the Computer Science Department, Technion--Israel Institute of Technology, Haifa 32000, Israel.}}


\maketitle

\begin{abstract}
A distinct difference configuration is a set of points in
$\mathbb{Z}^2$ with the property that the vectors (\emph{difference
vectors}) connecting any two of the points are all distinct. Many
specific examples of these configurations have been previously
studied: the class of distinct difference configurations includes both
Costas arrays and sonar sequences, for example.

Motivated by an application of these structures in key predistribution
for wireless sensor networks, we define the $k$-hop coverage of a
distinct difference configuration to be the number of distinct vectors
that can be expressed as the sum of $k$ or fewer difference
vectors. This is an important parameter when distinct difference
configurations are used in the wireless sensor application, as this
parameter describes the density of nodes that can be reached by a
short secure path in the network. We provide upper and lower bounds
for the $k$-hop coverage of a distinct difference configuration with
$m$ points, and exploit a connection with $B_{h}$ sequences to
construct configurations with maximal $k$-hop coverage. We also
construct distinct difference configurations that enable all small
vectors to be expressed as the sum of two of the difference vectors of
the configuration, an important task for local secure connectivity in
the application.
\end{abstract}

\section{Introduction}
\label{sec:introduction}

\IEEEPARstart{A}{} {\em distinct difference configuration}
$\DD(m)$ is a set of $m$ dots in a square grid, with the property
that the lines joining distinct pairs of dots are all different in
length or slope. For instance, the dots depicted in the following
array form a $\DD(3)$: {\tiny\begin{equation*}
\setlength{\arraycolsep}{.4\arraycolsep}
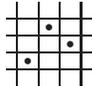
\begin{array}{c|c|c|c|c}
&&&&\\ \hline
&&\bullet&\\ \hline
 &&&\bullet\\ \hline
&\bullet &&&\\ \hline
&&&&
\end{array}
\end{equation*}}
If we pick a position on the square grid to be the origin, we may
think of the dots in a $\DD(m)$ as a set
$\{\mathbf{v}_1,\mathbf{v}_2,\dotsc,\mathbf{v}_m\}$ of vectors in
$\bZ^2$. The condition that the dots form a $\DD(m)$ is then the same
as the condition that the \emph{difference vectors}
$\mathbf{v}_i-\mathbf{v}_j$ with $i\neq j$ are all distinct. So we may
think of the dots in the example above as the set
$\{(0,0),(1,2),(2,1)\}$ of vectors; it is easy to verify that
the six difference vectors are all different in this case.

Many special classes of distinct difference configurations have been
studied previously: these include $B_2$ sequences over $\bZ$ and
Golomb rulers in the one-dimensional case, and Costas arrays, Golomb
rectangles and sonar sequences in the two-dimensional
case. See~\cite{BEMP1} for a summary of these configurations.

This paper is concerned with the \emph{$k$-hop properties} of
distinct difference configurations. Before we explain this, we
first need to discuss an application to key predistribution in
grid-based wireless sensor networks due to Blackburn, Etzion,
Martin and Paterson \cite{BEMP} that motivates our work.
\IEEEpubidadjcol
\subsection{Wireless Sensor Networks}

A {\em wireless sensor network} is a large collection of small sensor nodes that
are equipped with wireless communication capability. Sensor nodes have limited communication range and thus data transmitted over the network is typically passed from node to node in a series of {\em hops} in order to reach its end destination. Such networks can be employed for a wide range of applications \cite{romer-design}, whether scientific, commercial, humanitarian or military.  The data being transmitted over the wireless medium is
frequently valuable or sensitive; hence, there is a need for cryptographic techniques to provide
data integrity, confidentiality and authentication.

On deployment, the sensor nodes aim to form a secure and connected network. In other words, we desire a significant proportion of nodes within communication range to share cryptographic keys. The nodes'
size limits their computational power and battery capacity, so it is
assumed that the sensor nodes are unable to use public key cryptography to
establish shared keys.  So symmetric cryptographic keys are preloaded onto each
node before deployment: methods for deciding which keys are assigned
to a node are known as key predistribution schemes
(see~\cite{cysurvey,keithframework,newsurvey} for surveys of this
subject). The sensor nodes are assumed to be highly vulnerable to compromise,
so a single key should not be given to too many nodes. A balancing
constraint is that each node can only store a limited number of
keys. The aim is to design an efficient and secure key predistribution
scheme so that a sensor node can establish secure wireless links with many
of its neighbours: it is important to establish as many short secure
links in the network as possible, since the nodes' capacity to relay
information is very limited.

Key predistribution schemes for wireless sensor networks generally assume that the precise location of nodes is not known before deployment, hence schemes such as \cite{eschenauer02} aim to provide reasonable levels of ``average'' connectivity across the entire network. However in many applications the location of sensor nodes can be determined prior to deployment. In such cases this knowledge can be used to improve the efficiency of the underlying key predistribution scheme. One such scenario is that of networks consisting of a
large number of sensor nodes arranged in a square grid. Grid-based networks can arise in many applications, including soil moisture sensing \cite{soilnet}, monitoring conditions in an orchard
\cite{nectarine}, and measuring the efficiency of water use during irrigation
\cite{irrigate}.

\subsection{Key Predistribution for a Grid-based Network}

In~\cite{BEMP} a key predistribution scheme for a grid-based network was proposed and analysed. This scheme was shown to be significantly more efficient than using general approaches such as that of \cite{eschenauer02}. We now discuss this scheme in more detail.

Although the
number of sensor nodes is evidently finite in practice, it is convenient to
model the physical location of the nodes by the set of points of
$\mathbb{Z}^2$.  The scheme in \cite{BEMP}  employs a distinct difference configuration to create a key
predistribution scheme in the following way.
\begin{scheme}\label{BEMPscheme}
Let $D=\{\mathbf{v}_1,\mathbf{v}_2,\dotsc,\mathbf{v}_m\}$ be a
distinct difference configuration. Allocate keys to nodes as follows:
\begin{itemize}
\item Label each node with its position in $\mathbb{Z}^2$.
\item For every `shift' $\mathbf{u}\in\bZ^2$, generate a key
$k_\mathbf{u}$ and assign $k_\mathbf{u}$ to the nodes labelled by
$\mathbf{u}+\mathbf{v}_i$, for\, $i=1,2,\dotsc,m$.
\end{itemize}
\end{scheme}
\IEEEpubidadjcol More informally, we can think of the scheme as
covering $\bZ^2$ with all possible translations of the dots in $D$. We
generate one key per translation, and assign that key to all dots in
the corresponding translation of $D$. Distributing keys in this manner
ensures that each node stores $m$ keys and each key is shared by $m$
nodes.  In addition, the distinct difference property of the
configuration implies that any pair of nodes shares at most one key,
since the vector representing the difference in two nodes' positions
can occur at most once as a difference vector of $D$.  This leads to
an efficient distribution of keys, since for a fixed number of stored
keys the number of distinct pairs of nodes that share a key is
maximised.

As an example, consider the distinct difference
configuration given at the start of this introduction. If we use this
configuration for key distribution in Scheme~\ref{BEMPscheme}, each
node stores three keys.  Figure~\ref{costasex} illustrates this key
distribution: each square in the grid represents a node, and each
symbol contained in a square represents a key possessed by that node.
The central square stores keys marked by the letters $A$, $B$ and $C$;
two further nodes share each of these keys, which are marked in bold.
Letters in standard type represent keys used to connect the central
node to one of its neighbours via a two-hop path, other keys are
marked in grey.  Note that we have only illustrated some of the keys;
the pattern of key sharing extends in a similar manner throughout the
entire network. See~\cite{BEMP} for a comparison of how Scheme~\ref{BEMPscheme}
outperforms related key predistribution schemes in the literature.

\begin{figure}
\centering
\setlength{\unitlength}{.34mm}
\scriptsize
\begin{picture}(130,110)
\multiput(10,5)(15,0){8}{\line(0,1){115}}
\multiput(5,10)(0,15){8}{\line(1,0){115}}

\put(55,55){\makebox(5,5){$\mathbf A$}}
\put(75,95){\makebox(5,5){$\mathbf A$}}
\put(95,75){\makebox(5,5){$\mathbf A$}}

\put(60,65){\makebox(5,5){$\mathbf B$}}
\put(80,45){\makebox(5,5){$\mathbf B$}}
\put(40,25){\makebox(5,5){$\mathbf B$}}

\put(64.5,59){\makebox(5,5){$\mathbf C$}}
\put(45,78){\makebox(5,5){$\mathbf C$}}
\put(26.5,41){\makebox(5,5){$\mathbf C$}}

\put(71.5,56){\makebox(5,5){\textcolor{lgray}{$D$}}}
\put(89.5,93){\makebox(5,5){\textcolor{lgray}{$D$}}}
\put(109,75){\makebox(5,5){\textcolor{lgray}{$D$}}}

\put(75,63.5){\makebox(5,5){\textcolor{lgray}{$E$}}}
\put(94,45){\makebox(5,5){\textcolor{lgray}{$E$}}}
\put(56.5,26){\makebox(5,5){\textcolor{lgray}{$E$}}}

\put(90,63){\makebox(5,5){\textcolor{lgray}{$F$}}}
\put(110,45){\makebox(5,5){\textcolor{lgray}{$F$}}}
\put(71.5,26){\makebox(5,5){\textcolor{lgray}{$F$}}}

\put(45,63){\makebox(5,5){\textcolor{lgray}{$G$}}}
\put(64.5,45){\makebox(5,5){\textcolor{lgray}{$G$}}}
\put(26.5,26){\makebox(5,5){\textcolor{lgray}{$G$}}}

\put(60,78){\makebox(5,5){\textcolor{lgray}{$H$}}}
\put(79.5,60){\makebox(5,5){\textcolor{lgray}{$H$}}}
\put(41.5,41){\makebox(5,5){\textcolor{lgray}{$H$}}}

\put(60,93){\makebox(5,5){\textcolor{lgray}{$I$}}}
\put(80,75){\makebox(5,5){\textcolor{lgray}{$I$}}}
\put(41.5,56){\makebox(5,5){\textcolor{lgray}{$I$}}}

\put(64,75){\makebox(5,5){\textcolor{lgray}{$J$}}}
\put(45,93){\makebox(5,5){\textcolor{lgray}{$J$}}}
\put(26.5,56){\makebox(5,5){\textcolor{lgray}{$J$}}}

\put(56,41){\makebox(5,5){\textcolor{lgray}{$K$}}}
\put(75,78){\makebox(5,5){\textcolor{lgray}{$K$}}}
\put(94,60){\makebox(5,5){\textcolor{lgray}{$K$}}}

\put(49,60){\makebox(5,5){\textcolor{lgray}{$L$}}}
\put(30,78){\makebox(5,5){\textcolor{lgray}{$L$}}}
\put(11.5,41){\makebox(5,5){\textcolor{lgray}{$L$}}}

\put(56.5,71){\makebox(5,5){\textcolor{lgray}{$M$}}}
\put(75,108){\makebox(5,5){\textcolor{lgray}{$M$}}}
\put(94.5,90){\makebox(5,5){\textcolor{lgray}{$M$}}}

\put(59.5,48){\makebox(5,5){\textcolor{lgray}{$N$}}}
\put(79,30){\makebox(5,5){\textcolor{lgray}{$N$}}}
\put(41.5,11){\makebox(5,5){\textcolor{lgray}{$N$}}}

\put(49,75){\makebox(5,5){$O$}}
\put(30,93){\makebox(5,5){$O$}}
\put(11.5,56){\makebox(5,5){$O$}}

\put(49,90){\makebox(5,5){\textcolor{lgray}{$P$}}}
\put(30,108){\makebox(5,5){\textcolor{lgray}{$P$}}}
\put(11.5,71){\makebox(5,5){\textcolor{lgray}{$P$}}}

\put(64,90){\makebox(5,5){\textcolor{lgray}{$Q$}}}
\put(45,108){\makebox(5,5){\textcolor{lgray}{$Q$}}}
\put(26.5,71){\makebox(5,5){\textcolor{lgray}{$Q$}}}

\put(79,90){\makebox(5,5){$R$}}
\put(60,108){\makebox(5,5){$R$}}
\put(41.5,71){\makebox(5,5){$R$}}

\put(109,90){\makebox(5,5){\textcolor{lgray}{$S$}}}
\put(90,108){\makebox(5,5){\textcolor{lgray}{$S$}}}
\put(71.5,71){\makebox(5,5){\textcolor{lgray}{$S$}}}

\put(105,108){\makebox(5,5){\textcolor{lgray}{$T$}}}
\put(86.5,71){\makebox(5,5){\textcolor{lgray}{$T$}}}

\put(109,105){\makebox(5,5){\textcolor{lgray}{$U$}}}
\put(71.5,86){\makebox(5,5){\textcolor{lgray}{$U$}}}

\put(71.5,41){\makebox(6,6){$V$}}
\put(89.5,78){\makebox(6,6){$V$}}
\put(109,60){\makebox(6,6){$V$}}

\put(74.5,48){\makebox(5,5){$W$}}
\put(94,30){\makebox(5,5){$W$}}
\put(57,11){\makebox(5,5){$W$}}

\put(45,48){\makebox(5,5){\textcolor{lgray}{$X$}}}
\put(64.5,30){\makebox(5,5){\textcolor{lgray}{$X$}}}
\put(26.5,11){\makebox(5,5){\textcolor{lgray}{$X$}}}

\put(45,33){\makebox(5,5){$Y$}}
\put(64.5,15){\makebox(5,5){$Y$}}

\put(49,45){\makebox(5,5){\textcolor{lgray}{$Z$}}}
\put(30,63){\makebox(5,5){\textcolor{lgray}{$Z$}}}
\put(11.5,26){\makebox(5,5){\textcolor{lgray}{$Z$}}}

\put(49,30){\makebox(5,5){$\Delta$}}
\put(30,48){\makebox(5,5){$\Delta$}}
\put(11.5,11){\makebox(5,5){\textcolor{lgray}{$\Delta$}}}

\put(15,63){\makebox(5,5){$\Phi$}}
\put(35.5,45){\makebox(5,5){$\Phi$}}

\put(49,15){\makebox(5,5){\textcolor{lgray}{$\Gamma$}}}
\put(30,33){\makebox(5,5){\textcolor{lgray}{$\Gamma$}}}

\put(60,33){\makebox(5,5){\textcolor{lgray}{$\Theta$}}}
\put(79.5,15){\makebox(5,5){\textcolor{lgray}{$\Theta$}}}

\put(105,63){\makebox(5,5){\textcolor{lgray}{$\Lambda$}}}
\put(86.5,26){\makebox(5,5){\textcolor{lgray}{$\Lambda$}}}

\put(105,78){\makebox(5,5){\textcolor{lgray}{$\Xi$}}}
\put(86.5,41){\makebox(5,5){\textcolor{lgray}{$\Xi$}}}

\put(89.5,48){\makebox(5,5){\textcolor{lgray}{$\Pi$}}}
\put(109,30){\makebox(5,5){\textcolor{lgray}{$\Pi$}}}
\put(72,11){\makebox(5,5){\textcolor{lgray}{$\Pi$}}}

\put(34.5,60){\makebox(5,5){\textcolor{lgray}{$\Sigma$}}}
\put(15,78){\makebox(5,5){\textcolor{lgray}{$\Sigma$}}}

\put(34.5,30){\makebox(5,5){\textcolor{lgray}{$\Upsilon$}}}
\put(15,48){\makebox(5,5){\textcolor{lgray}{$\Upsilon$}}}

\put(65,105){\makebox(5,5){\textcolor{lgray}{$\Psi$}}}
\put(26.5,86){\makebox(5,5){\textcolor{lgray}{$\Psi$}}}

\put(80,105){\makebox(5,5){\textcolor{lgray}{$\Omega$}}}
\put(41.5,86){\makebox(5,5){\textcolor{lgray}{$\Omega$}}}

\put(95,105){\makebox(5,5){\textcolor{lgray}{$\digamma$}}}
\put(56.5,86){\makebox(5,5){\textcolor{lgray}{$\digamma$}}}

\thicklines
\put(55,55){\framebox(14.8,14.8){}}
\end{picture}
\caption{Key distribution using a distinct difference
configuration.}\label{costasex}
\end{figure}

Note that the sensors' strictly limited battery power limits the range over
which they can feasibly communicate. In support of Scheme~\ref{BEMPscheme},
distinct difference configurations with bounds on the distance between
any two dots in the configuration were considered in \cite{BEMP}. Supposing that each sensor has a
fixed communication range $r$, a {\em $\DD(m,r)$} is defined to be a
$\DD(m)$ in which the Euclidean distance between any two points of
the configuration is at most $r$. From an application point of view,
it is only necessary for a pair of nodes to share a key if they are
located within communication range of each other; the use of a
$\DD(m,r)$ in Scheme~\ref{BEMPscheme} ensures that this is the case.

While Scheme~\ref{BEMPscheme} was designed to suit wireless sensor
networks in which the sensors are arranged in a square grid, for
certain applications a hexagonal arrangement of sensor nodes may be
preferred, as it yields the most efficient packing of sensors (see
\cite{CoSl93} for details of circle packings in the
plane). Section~\ref{sec:models} defines the hexagonal model more
precisely and discusses the relationship between the two models.
Scheme~\ref{BEMPscheme} is easily adapted to suit sensors arranged in
a hexagonal grid by replacing the $\DD(m)$ by a $\DD^*(m)$, which we
informally define to be a set of $m$ dots on a hexagonal grid such
that the vector differences between pairs of dots are distinct.  We
define a $\DD^*(m,r)$ to be a $\DD^*(m)$ in which the Euclidean
distance between any pair of dots is at most~$r$.  Another model
that is natural when working with either the square or hexagonal grids
is to replace the Euclidean metric by its discrete equivalent: the
Manhattan metric (in the case of square grids), or an analogous metric
on the hexagonal grid; in this case, we use the notation $\oDD(m,r)$
and $\oDD^*(m,r)$, respectively.  Constructions and bounds on the
parameters for such configurations were studied in \cite{BEMP1}.
Section~\ref{sec:models} contains a summary of the relationships
between configurations based on different grids when using different
metrics.

\subsection{Contributions}

Recall that wireless sensor networks rely on data being relayed via intermediate nodes using a series of hops. From an efficiency perspective it is thus of interest to consider properties relating to the nodes that can be reached from a specific node by means of a restricted number of hops.

If two nodes $A$ and $B$
are within communication range and share a key we say there is a {\em
one-hop path} between $A$ and $B$.  If they do not share a key,
however, they may still be able to establish a secure connection if
there is a node $C$ that is within range of $A$ and $B$ and shares a
key with each of them.  This is referred to as a {\em two-hop path};
more generally we consider {\em $k$-hop paths} of the form
$A-C_1-C_2\dotsc-C_{k-1}-B$, where there is a one-hop path between any
two adjacent users in the chain.  A significant, and widely studied,
measure of the performance of a key predistribution scheme for a
wireless sensor network is the expected number of nodes with which a
given node can communicate via a one hop or two-hop path (we do not
count the given node in this total).  As in \cite{BEMP}, we refer to
this parameter as the {\em two-hop coverage} of the scheme.  More
generally, we can define the {\em $k$-hop coverage} to be the expected
number of nodes with which a given node can communicate via some
$\ell$-hop path with $1\leq \ell\leq k$ (where we do not count the
given node itself).

This measure is important from the point of view
of our application, since it captures the ability of the network to
transmit information in the context of the nodes' limited capacity to
relay messages. The case when $k=2$ is the most studied situation in
the literature, since results are often easier to establish than in the
general $k$-hop case.  Lee and Stinson use the notation ${\rm \bf
Pr}_1+{\rm \bf Pr}_2$ to describe this quantity, referring to it as
the {\em local connectivity} \cite{leestin07}; similar metrics are used in
\cite{Du03apairwise,Chan03randomkey}, and various related measures of
the expected number of hops required for secure communication between
two nodes are prevalent in the sensor network literature
\cite{eschenauer02, expander, liuningli05}.

We define the $k$-hop coverage of a distinct difference configuration to be the $k$-hop
coverage of the resulting instance of Scheme~\ref{BEMPscheme}. In \cite{BEMP} a number of distinct difference configurations with good two-hop coverage were found by computer search. However no concrete construction techniques were provided. In this paper we provide an exposition of the two-hop coverage case, as well as consider the generalisation to $k$-hop coverage.

Section~\ref{sec:twohop}
is devoted to a study of the $k$-hop coverage
$C_k(D)$ obtained by the use of the distinct difference configuration
$D=\{\mathbf{v}_1,\mathbf{v}_2,\ldots,\mathbf{v}_m\}$ in
Scheme~\ref{BEMPscheme}.  Subsection~\ref{sec:khop_characterization}
shows how to calculate the $k$-hop coverage from the vectors
$\mathbf{v}_1,\mathbf{v}_2,\ldots,\mathbf{v}_m$. In
Subsection~\ref{sec:max} we study configurations where $C_k(D)$
is as large as possible, and show a connection between such
configurations and $B_h$ sequences (a well studied concept in
combinatorial number theory). We determine the maximum value of the
$k$-hop coverage $C_k(D)$ where $D$ is a $\DD(m)$ (or a $\DD^*(m)$), and
show that $D$ achieves this level of $k$-hop coverage if and only if
$D$ is a $B_{2k}$ sequence. If we restrict $D$ to be a $\DD(m,r)$ for
some small integer $r$, we might no longer be able to achieve this
maximum value of $C_k(D)$: we provide bounds on the smallest value of
$r$ for which there exists a configuration $D$ which is a $\DD(m,r)$
with $C_k(D)$ maximal. We also provide similar bounds on this smallest
value of $r$ when we consider configurations $\DD^*(m,r)$ in the
hexagonal grid. Finally, in Subsection~\ref{sec:min}, we provide a
lower bound on $C_k(D)$ and characterise those configurations that
meet this lower bound.

Using a distinct difference configuration with maximal $k$-hop
coverage ensures that as many users as possible are connected by
$k$-hop paths.  However, in many applications these paths are used to
establish keys which are later used for direct communication between
the two end nodes: thus we are only interested in $k$-hop paths whose
start and end nodes are within communication range. For these
applications, rather than optimising the total number of pairs of
users connected by $k$-hop paths we wish to optimise coverage in a
locally defined region: We say that a $\DD(m)$ or $\DD^*(m)$ achieves
{\em complete $k$-hop coverage with respect to a region $R$ and point
$\mathbf{p}\in R$} if every point in $R$ can be reached by a two-hop
path from~$\mathbf{p}$. This means that every node $\mathbf{u}$ can
communicate via a $k$-hop path with the nodes in the region
corresponding to a shift of $R$ that moves $\mathbf{p}$ to
$\mathbf{u}$, giving Scheme~\ref{BEMPscheme} good local
connectivity. In Section~\ref{sec:complete-two-hop} we give a
construction for a $\DD(m)$ that achieves complete two-hop coverage
with respect to the centre of a $(2p-3)\times(2p-1)$ rectangle when
$p$ is prime.

\section{Different Grids and Different Metrics}
\label{sec:models}
\subsection{Square and Hexagonal Grids}

Suppose that the sensor nodes are arranged in a square grid, and the
shortest distance between a pair of nodes is $1$. So we tile the plane
by unit squares, and think of the nodes as lying at the centres of
these squares. By supposing one of the nodes is at the origin, the
location of a node can be identified with a vector in $\bZ^2$. Because
of this, we call $\bZ^2$ the \emph{square grid}.

A hexagonal arrangement of sensor nodes is obtained by tiling the
plane with regular hexagons and placing a node at the centre of each
hexagon. We suppose that one of the nodes is located at the origin and
the shortest distance between two nodes is 1. In a similar way to the
square grid, the locations of the nodes can be represented by vectors
in the set $\Lambda_H=\{\lambda(1,0)+\mu(-1/2,\sqrt{3}/2)\vert
\lambda,\mu\in\mathbb{Z}\}$, which we call the \emph{hexagonal
grid}.

We have already defined a (square) distinct difference configuration
$\DD(m)$ to be a set $D=\{\mathbf{v}_1,
\mathbf{v}_2,\ldots,\mathbf{v}_m\}\subseteq \bZ^2$ of $m$ dots with the
property that the difference vectors $\mathbf{v}_i-\mathbf{v}_j$ for
$i\not=j$ between any pair of dots are distinct. In the same way, we
define a (hexagonal) distinct difference configuration $\DD^*(m)$ to
be a set $D=\{\mathbf{v}_1,
\mathbf{v}_2,\ldots,\mathbf{v}_m\}\subseteq \Lambda_H$ of $m$ dots in the
hexagonal grid with the property that the difference vectors
$\mathbf{v}_i-\mathbf{v}_j$ for $i\not=j$ are distinct. A hexagonal
distinct difference configuration can be used in
Scheme~\ref{BEMPscheme} for sensors arranged in a hexagonal grid,
provided that shifts $\mathbf{u}\in\Lambda_H$ are used: as in the
square grid, every node is assigned $m$ keys and the distinct
difference property implies that any pair of nodes has at most one key
in common. We define a $\DD^*(m,r)$ to be a $\DD^*(m)$ in which the
Euclidean distance between any pair of dots in the configuration is at
most $r$: these configurations must be used when the wireless
communication range of a sensor node is $r$.

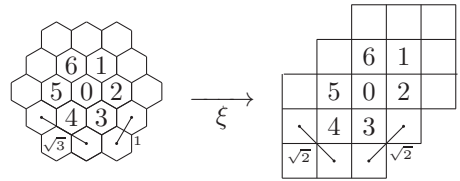
\begin{figure}[t]
\centering
\setlength{\unitlength}{.4mm}
\begin{picture}(150,70)(-10,-40)
\multiput(10,5)(10,0){3}{
\put(-10,0){\line(0,1){5.77350269}}
\put(0,0){\line(0,1){5.77350269}}
\put(0,5.77350269){\line(-1732,1000){5}}
\put(-10,5.77350269){\line(1732,1000){5}}
}
\multiput(5,-3.66033872)(10,0){4}{
\put(-10,0){\line(0,1){5.77350269}}
\put(0,0){\line(0,1){5.77350269}}
\put(0,5.77350269){\line(-1732,1000){5}}
\put(-10,5.77350269){\line(1732,1000){5}}
}
\multiput(0,-12.3206774)(10,0){5}{
\put(-10,0){\line(0,1){5.77350269}}
\put(0,0){\line(0,1){5.77350269}}
\put(0,5.77350269){\line(-1732,1000){5}}
\put(-10,5.77350269){\line(1732,1000){5}}
}
\multiput(5, -20.9810162)(10,0){4}{
\put(-10,0){\line(0,1){5.77350269}}
\put(0,0){\line(0,1){5.77350269}}
\put(0,5.77350269){\line(-1732,1000){5}}
\put(-10,5.77350269){\line(-1732,1000){5}}
\put(-10,5.77350269){\line(1732,1000){5}}
\put(0,5.77350269){\line(1732,1000){5}}
}
\multiput(10, -29.6413549)(10,0){3}{
\put(-10,0){\line(0,1){5.77350269}}
\put(0,0){\line(0,1){5.77350269}}
\put(0,5.77350269){\line(-1732,1000){5}}
\put(-10,5.77350269){\line(-1732,1000){5}}
\put(-10,5.77350269){\line(1732,1000){5}}
\put(0,5.77350269){\line(1732,1000){5}}
}
\multiput(15, -38.3016936)(10,0){2}{
\put(0,5.77350269){\line(-1732,1000){5}}
\put(-10,5.77350269){\line(-1732,1000){5}}
\put(-10,5.77350269){\line(1732,1000){5}}
\put(0,5.77350269){\line(1732,1000){5}}}

\put(15,-9.43392605){\makebox(0,0){$0$}}
\put(5,-9.43392605){\makebox(0,0){$5$}}
\put(25,-9.43392605){\makebox(0,0){$2$}}
\put(10,-0.77358733){\makebox(0,0){$6$}}
\put(20,-0.77358733){\makebox(0,0){$1$}}
\put(10,-18.0942648){\makebox(0,0){$4$}}
\put(20,-18.0942648){\makebox(0,0){$3$}}
\put(0,-18.0942648){\line(172,-100){15}}
\put(0,-18.0942648){\circle*{1}}
\put(30,-18.0942648){\circle*{1}}
\put(15,-26.7546035){\circle*{1}}
\put(25,-26.7546035){\circle*{1}}
\put(30,-18.0942648){\line(-100,-172){5}}

\put(85.7735027,-20.9809314){\line(1,-1){11.5470054}}
\put(120.414519,-20.9809314){\line(-1,-1){11.5470054}}
\put(85.7735027,-20.9809314){\circle*{1}}
\put(120.414519,-20.9809314){\circle*{1}}
\put(97.3205081,-32.5279368){\circle*{1}}
\put(108.867514,-32.5279368){\circle*{1}}
\thinlines

\put(4,-27){\makebox(0,0){\tiny$\sqrt{3}$}}
\put(32,-26){\makebox(0,0){\tiny$1$}}
\put(85.3205081,-30.9809314){\makebox(0,0){\tiny$\sqrt{2}$}}
\put(119.3205081,-29.9809314){\makebox(0,0){\tiny$\sqrt{2}$}}

\put(60,-16.43392605){\makebox(0,0){$\overrightarrow{\quad\xi\quad}$}}

\put(80,-38.3014395){\line(1,0){34.6410162}}
\put(80,-26.7544341){\line(1,0){46.1880215}}
\put(80,-15.2074287){\line(1,0){57.7350269}}
\put(80.5470054,-3.66042332){\line(1,0){57.7350269}}
\put(91.5470054,7.88658206){\line(1,0){46.1880215}}
\put(103.094011, 19.4335874){\line(1,0){34.6410162}}

\put(80,-38.3014395){\line(0,1){34.6410162}}
\put(91.5470054,-38.3014395){\line(0,1){46.1880215}}
\put(103.094011,-38.3014395){\line(0,1){57.7350269}}
\put(114.641016,-38.3014395){\line(0,1){57.7350269}}
\put(126.188022,-26.7544341){\line(0,1){46.1880215}}
\put(137.735027,-15.2074287){\line(0,1){34.6410162}}

\put(108.867513,-9.43392605){\makebox(0,0){$0$}}
\put( 97.3205081,-9.43392605){\makebox(0,0){$5$}}
\put(120.414519,-9.43392605){\makebox(0,0){$2$}}
\put(108.867513,2.11307933){\makebox(0,0){$6$}}
\put(97.3205081,-20.9809314){\makebox(0,0){$4$}}
\put(120.414519,2.11307933){\makebox(0,0){$1$}}
\put(108.867513,-20.9809314){\makebox(0,0){$3$}}
\end{picture}
\caption{A transformation from a hexagonal grid to a square grid (grid points are represented by the centres of the cells).} \label{fig:hexmodel}
\end{figure}

The map $\xi\colon \mathbb{R}^2 \rightarrow \mathbb{R}^2$
defined by
\begin{equation*}
\xi\colon(x,y)\mapsto (x+\frac{y}{\sqrt{3}},\frac{2y}{\sqrt{3}})
\end{equation*}
induces a bijection from $\Lambda_H$ to $\mathbb{Z}$. This is
illustrated in Fig.~\ref{fig:hexmodel}, in which the cells whose
centres form the points of the grid are depicted.
We can use $\xi$ and $\xi^{-1}$ to convert a $\DD^*(m)$ into a
$\DD(m)$ and {\it vice versa}:
\begin{theorem}
\label{thm:DDequivalence}
If $D=\{\mathbf{v}_1,\mathbf{v}_2,\dotsc,\mathbf{v}_m\}$ is a
$\DD^*(m)$, then
$\xi(D)=\{\xi(\mathbf{v}_1),\xi(\mathbf{v}_2),\dotsc,\xi(\mathbf{v}_m)\}$
is a $\DD(m)$.  Similarly, if $D^\prime$ is a $\DD(m)$, then
$\xi^{-1}(D^\prime)$ is a $\DD^*(m)$.
\end{theorem}
\begin{IEEEproof}
Since $\xi$ is a linear bijection, we have that
$\mathbf{v}_i-\mathbf{v}_j=\mathbf{v}_k-\mathbf{v}_\ell$ if and only
if
$\xi(\mathbf{v}_i)-\xi(\mathbf{v}_j)=\xi(\mathbf{v}_k)-\xi(\mathbf{v}_\ell)$;
the first statement of the theorem follows directly.  The second
statement follows as $\xi^{-1}$ is also a linear bijection.
\end{IEEEproof}

Despite Theorem~\ref{thm:DDequivalence}, the square and hexagonal
models differ once we are interested in distances between dots, since $\xi$
does not preserve Euclidean distances. Fig.~\ref{fig:hexmodel} shows a
line segment of length $\sqrt{3}$ that transforms into one of length
$\sqrt{2}$, and one of length $1$ that also transforms into one of
length $\sqrt{2}$.  It is straightforward to show that these line
segments represent the maximum extent to which $\xi$ can extend or
contract the length of a vector; we formalise this in the following
theorem:
\begin{theorem}\label{thm:hextosquare}
If $D$ is a $\DD^*(m,r)$ then $\xi(D)$ is a $\DD(m,r\sqrt{2})$.  If
$D^\prime$ is a $\DD(m,r)$, then $\xi^{-1}(D^\prime)$ is a
$\DD^*(m,r\sqrt{3/2})$.
\end{theorem}
Thus we can convert between results about $\DD(m,r)$ and results about
$\DD^*(m,r)$ (although the bounds on the converted lengths are not
tight in general).
\subsection{Alternative Metrics on Grids}
In \cite{BEMP}, the need to take sensor nodes' communication range
into account when using distinct difference configurations to
distribute keys to sensors arranged in a square grid motivated the
definition of a $\DD(m,r)$ based on a Euclidean measure of
distance.  However, when working with a square grid it is natural to
consider the Manhattan metric (also known as the Lee metric), in which
the distance between dots with coordinates $(i_1,j_1)$ and
$(i_2,j_2)$ is given by $\lvert i_1-i_2\rvert+\lvert
j_1-j_2\rvert$. Distinct difference configurations
$\oDD(m,r)$ in which the distance between dots in the
configuration is at most $r$ in the Manhattan metric were studied in
\cite{BEMP1}.  A ball of radius $r$ in this metric is referred to as a
{\em Lee sphere} (Fig.~\ref{leeee}), and for small $r$ gives a
reasonable approximation of a Euclidean circle.  The well-known
relation between these two metrics is expressed in the following
theorem, which permits conversion between results about
$\oDD(m,r)$ and results about $\DD(m,r)$.
\begin{figure}[tb]
\centering
\subfloat[][Lee sphere of radius 2]{
\begin{picture}(90,50)(-77,0)
\put(-36,5){\framebox(8,40){}}

\put(-52,21){\framebox(40,8){}}

\put(-44,13){\framebox(24,24){}}
\end{picture}
\label{leeee}}
\subfloat[][Hexagonal ball of radius 1]{
\setlength{\unitlength}{.4mm}
\begin{picture}(80,40)(55,0)
\multiput(86.3333333,30.980762)(13.3333333,0){3}{\line(0,1){7.69800359}}

\multiput(93,42.5277674)(13.3333333,0){2}{\line(1732,-1000){6.6666667}}

\multiput(93,42.5277674)(13.3333333,0){2}{\line(-1732,-1000){6.6666667}}

\multiput(79.6666666,19.4337566)(13.3333333,0){4}{\line(0,1){7.69800359}}

\multiput(86.3333333,30.980762)(13.3333333,0){3}{\line(1732,-1000){6.6666667}}

\multiput(86.3333333,30.980762)(13.3333333,0){3}{\line(-1732,-1000){6.6666667}}

\multiput(86.3333333,7.88675122)(13.3333333,0){3}{\line(0,1){7.69800359}}

\multiput(79.6666667,19.4337566)(13.3333333,0){3}{\line(1732,-1000){6.6666667}}

\multiput(93,19.4337566)(13.3333333,0){3}{\line(-1732,-1000){6.6666667}}

\multiput(86.3333333,7.88675122)(13.3333333,0){2}{\line(1732,-1000){6.6666667}}

\multiput(99.666667,7.88675122)(13.3333333,0){2}{\line(-1732,-1000){6.6666667}}

\end{picture}
}
\end{figure}
\begin{theorem}\label{euclee}
For $r\in\mathbb{Z}$, a $\oDD(m,r)$ is a $\DD(m,r)$ and a $\DD(m,r)$
is a $\oDD(m,\lceil\sqrt{2}r\rceil)$.
\end{theorem}

For the hexagonal grid, we say that a given point is {\em adjacent} to
the six grid points that lie at Euclidean distance 1 from that point
(for example, in Fig.~\ref{fig:hexmodel} the points at the centres of
cells $1,2,\dotsc,6$ are adjacent to the point at the centre of cell
0).  We can then define a graph in which the grid points correspond to
vertices, with edges connecting vertices whose grid points are
adjacent.  This gives rise to a {\em hexagonal metric} in which the
distance between two points is the length of the shortest path between
the corresponding vertices in the graph.  A distinct difference
configuration in which the hexagonal distance between any two points
is at most $r$ is denoted $\oDD^*(m,r)$.  The relation between the
hexagonal and Euclidean metrics can be used to prove the following
theorem:
\begin{theorem}
\label{thm:euchex}
For $r\in\mathbb{Z}$, a $\oDD^*(m,r)$ is a $\DD^*(m,r)$ and a
$\DD^*(m,r)$ is a $\oDD^*(m,\lceil\frac{2}{\sqrt{3}}r\rceil)$.
\end{theorem}
We note that the hexagonal metric gives a closer approximation to the
Euclidean distance than the Manhattan metric.

\section{$k$-Hop Coverage}
\label{sec:twohop}

In this section we investigate the properties of distinct difference configurations with respect to their $k$-hop coverage. While the motivation for this work comes from the application, the results are of independent combinatorial interest.

\subsection{Characterising $k$-hop coverage}
\label{sec:khop_characterization}

Let $D$ be a (square or hexagonal) distinct difference configuration
given by
$D=\{\mathbf{v}_1,\mathbf{v}_2,\dotsc,\mathbf{v}_m\}$. Define
$C_k(D)$ to be the number of non-zero vectors that can be written as
the sum of $k$ or fewer difference vectors. So $C_k(D)$ is the number
of non-zero vectors of the form
\begin{equation}
\label{eqn:CkD}
\sum\limits_{i=1}^\ell(\mathbf{v}_{\alpha_i}-\mathbf{v}_{\beta_i})
\end{equation}
where $\alpha_i,\beta_i \in \{1,2,\dotsc,m\}$ with
$\alpha_i\not=\beta_i$ and where $0\leq\ell\leq k$.
\begin{theorem}
\label{thm:khop_character}
Suppose that $D$ is used in Scheme~\ref{BEMPscheme}. Then the $k$-hop
coverage of the scheme is equal to $C_k(D)$.
\end{theorem}
\begin{IEEEproof}
Let $\mathbf{x}$ be any fixed node. Two nodes that share a key are
located at points of the form $\mathbf{v}_i+\mathbf{u}$ and
$\mathbf{v}_j+\mathbf{u}$ for some $i,j\in\{1,2,\dotsc,m\}$ and some
shift $\mathbf{u}$.  This implies that the vector difference between
their positions is $\mathbf{v}_i-\mathbf{v}_j$, which is a difference
vector of $D$. Hence a one-hop path between nodes with keys
distributed according to Scheme~\ref{BEMPscheme} corresponds to a
difference vector of the underlying distinct difference configuration.
So there is an $\ell$-hop path from $\mathbf{x}$ to another node
$\mathbf{y}$ if and only if the vector difference between their
positions is the sum of $\ell$ difference vectors.  Note also that
$\mathbf{x}=\mathbf{y}$ if and only if this sum is the zero vector:
since we do not count $\mathbf{x}$ in the $k$-hop coverage, we are
only interested in sums of the form~\eqref{eqn:CkD} which are
non-zero. So $C_k(D)$ is equal to the $k$-hop coverage of
Scheme~\ref{BEMPscheme} implemented using $D$, as required.
\end{IEEEproof}

\begin{theorem}\label{thm:khopsame}
Let $\xi:\bR^2\rightarrow\bR^2$ be the map defined in
Section~\ref{sec:models}. Let $D$ be a $\DD^*(m)$ and let $D'$ be a
$\DD(m)$ such that $D'=\xi(D)$. Then the $k$-hop coverage of $D$ is
equal to the $k$-hop coverage of $D'$.
\end{theorem}
\begin{IEEEproof}
Theorem~\ref{thm:khop_character} shows that we must show that
$C_k(D)=C_k(D')$. But $C_k(D)$ and $C_k(D')$ both count the number of
non-zero vectors that can be expressed as the sum of $k$ or fewer
difference vectors (of $D$ or $D'$ respectively). The theorem now
follows, since $\xi$ is a linear bijection.
\end{IEEEproof}

\subsection{Maximal $k$-hop coverage}
\label{sec:max}

In this subsection we determine the maximal $k$-hop coverage of a
$\DD(m)$.  By Theorem~\ref{thm:khopsame}, these results apply equally
to a $\DD^*(m)$.  We begin with some preliminary notation and lemmas.

For a non-negative integer $k$ we define a set $H_k$ of $m$-tuples of
integers as follows:
\begin{multline*}
H_k= \Bigg\{(a_1,a_2,\ldots,a_m)\in\mathbb{Z}^m\Bigg\vert\\ \sum_{i=1}^m a_i=0,\ \sum_{\{i:a_i>0\}}a_i=k\Bigg \}.
\end{multline*}
For example, when $m=3$ the triple $(0,0,0)$ is the unique element of
$H_0$, the triple $(1,-1,0)$ is a typical element of $H_1$, and the
triples $(2,-2,0)$, $(2,-1,-1)$ and $(1,1,-2)$ are typical elements of
$H_2$.  The following results about the sets $H_k$ are easily proved.
\begin{lemma}\label{lem:facts}
Define the sets $H_k$ as above.
\begin{enumerate}[(i)]
\item\label{fact1} Let $\mathbf{a}\in H_{k_1}$ and $\mathbf{b}\in H_{k_2}$. Then $\mathbf{a}+\mathbf{b}\in
H_{k_3}$ where $k_3$ is an integer satisfying $0\leq k_3\leq
k_1+k_2$. In particular, if a non-zero $m$-tuple $\mathbf{v}$ is a sum of $k$
$m$-tuples from $H_1$, then $\mathbf{v}\in H_{k_3}$ for some $k_3$ satisfying
$1\leq k_3\leq k$.
\item\label{fact2} Let $\mathbf{a}\in H_{k_1}$ and $\mathbf{b}\in H_{k_2}$ with $\mathbf{a}\neq\mathbf{b}$. Then
$\mathbf{a}-\mathbf{b}\in H_{k_3}$ where $k_3$ is an integer satisfying $1\leq k_3\leq
k_1+k_2$.
\item\label{fact3} Any element of $H_{k_1}$ may be written as the sum of
$k_1$ elements from $H_1$.
\end{enumerate}
\end{lemma}

The connection between $H_k$ and the $k$-hop coverage of $\DD(m)$ is given by the following theorem:

\begin{theorem}
\label{thm:khop_coverage} The $k$-hop coverage of a $\DD(m)$ is
at most $\sum_{i=1}^k |H_i|$, with equality if and only if all the
vectors $\sum_{i=1}^m a_i \mathbf{v}_i$ with $(a_1,a_2,\dotsc,a_m)\in
\bigcup_{j=0}^k H_j$ are distinct.
\end{theorem}
\begin{IEEEproof}
The difference vectors of $D$ are precisely the vectors of the form
$\sum_{i=1}^m a_i\mathbf{v}_i$ where $\mathbf{a}\in H_1$. By
Lemma~\ref{lem:facts}~(i) and~(iii), a vector is a sum of $k$ or fewer
difference vectors if and only if it can be written in the form
$\sum_{i=1}^m a_i \mathbf{v}_i$ with $(a_1,a_2,\dotsc,a_m)\in
\bigcup_{j=0}^k H_j$. The zero vector can always be written in this
form, since the sum is zero when $(a_1,a_2,\dotsc,a_m)\in H_0$. Since
and we are only interested in non-zero vectors, we find that
\begin{align*}
C_k(D)+1&= \left|\left\{\sum_{i=1}^m a_i \mathbf{v}_i\text{ where } \mathbf{a}\in \bigcup_{j=0}^k H_j\right\}\right|\\
&\leq\left(\sum_{i=0}^k |H_i|\right)\\
&=1+\left(\sum_{i=1}^k |H_i|\right),
\end{align*}
and it is clear that equality is satisfied if and only if the vectors
$\sum_{i=1}^m a_i \mathbf{v}_i\text{ where } \mathbf{a}\in
\bigcup_{j=0}^k H_j$ are distinct. Thus the theorem follows.
\end{IEEEproof}
\begin{corollary}
\label{cor:twohop_coverage} The two-hop coverage of a $\DD(m)$ is
at most
\begin{multline*}
\frac{1}{4}m(m-1)(m-2)(m-3)+m(m-1)(m-2)\\
+2m(m-1)=\frac{1}{4}m(m-1)(m^2-m+6).
\end{multline*}
\end{corollary}
\begin{IEEEproof}
By Theorem~\ref{thm:khop_coverage} the two-hop coverage is at most
$\lvert H_1\rvert+\lvert H_2\rvert$.  It is clear that $|H_1|=m(m-1)$,
since the $m$-tuples in $H_1$ have exactly two non-zero components,
one equal to $1$ and one equal to $-1$. To determine $|H_2|$, note
that there are four types of element in $H_2$, corresponding to the
four possibilities for the multiset of non-zero coefficients in an
$m$-tuple $\mathbf{a}\in H_2$ (see Table~\ref{countH2_table}). The
number of elements in $H_2$ of each type is equal to $(1/s)m!/(m-t)!$,
where $t$ is the number of non-zero components in an $m$-tuple of this
type, and $s$ is the number of symmetries that preserve such
$m$-tuples.
\begin{table}[t]
\renewcommand{\arraystretch}{1.4}
\caption{Counting elements in $H_2$} \label{countH2_table}
\centering
$\begin{array}{cccc}
\hline
\text{Type}&\text{Non-zero coeffs}&\text{Symm}&\text{Number}\\\hline
\text{(a)}&1,1,-1,-1&4&\frac{1}{4}m(m-1)(m-2)(m-3)\\
\text{(b)}&2,-1,-1&2&\frac{1}{2}m(m-1)(m-2)\\
\text{(c)}&1,1,-2&2&\frac{1}{2}m(m-1)(m-2)\\
\text{(d)}&2,-2&1&m(m-1)\\
\hline
\end{array}$
\end{table}
Thus $|H_2|=\frac{1}{4}m(m-1)(m-2)(m-3)+m(m-1)(m-2)+m(m-1)$, and
so the bound of the corollary follows.
\end{IEEEproof}

In order to show that the bound of Theorem~\ref{thm:khop_coverage} and
Corollary~\ref{cor:twohop_coverage} is tight, we must show that there
exists a $\DD(m)$ given by
\makebox{$D=\{\mathbf{v}_1,\mathbf{v}_2,\dotsc, \mathbf{v}_m\}$} such
that the vectors $\sum_{i=1}^ma_i\mathbf{v}_i$, where $\mathbf{a}\in
H_0\cup H_1 \cup \cdots \cup H_k$, are all distinct.  This is not
difficult to do: for example we may choose $\mathbf{v}_i=((2k+1)^i,0)$
for $i=1,2,\dotsc,m$. We say that a configuration meeting the bound of
Theorem~\ref{thm:khop_coverage} has {\em maximal $k$-hop
coverage}. Note that the example we have just given of a configuration
with maximal $k$-hop coverage is not useful for our application, as
the dots in the configuration are exponentially far apart: we would
like to construct a $\DD(m,r)$ with $r$ small having maximal $k$-hop
coverage. In order to do this, we now aim to characterise those
configurations with maximal $k$-hop coverage in terms of the much
studied concept of $B_h$ sequences (see below).  First, we make the
following observation.
\begin{lemma}
\label{lem:maximalkhop_criterion}
The $k$-hop coverage of a $\DD(m)$ given by
$D=\{\mathbf{v}_1,\mathbf{v}_2,\dotsc,\mathbf{v}_m\}$ meets the bound
of Theorem~\ref{thm:khop_coverage} if and only if
$\sum_{i=1}^mc_i\mathbf{v}_i\neq0$ for all $\mathbf{c}\in
\bigcup_{i=1}^{2k} H_i$.
\end{lemma}
\begin{IEEEproof}
Suppose that $D$ does not meet the bound of
Theorem~\ref{thm:khop_coverage}.  Then Theorem~\ref{thm:khop_coverage}
implies that
$\sum_{i=1}^ma_i\mathbf{v}_i=\sum_{i=1}^mb_i\mathbf{v}_i$, where
$\mathbf{a},\mathbf{b}\in \bigcup_{i=0}^k H_i$ and $\mathbf{a}\neq
\mathbf{b}$. Writing $\mathbf{c}=\mathbf{a}-\mathbf{b}$ we have
that $\sum_{i=1}^mc_i\mathbf{v}_i=0$, and $\mathbf{c}\in
\bigcup_{i=1}^{2k} H_i$ by Lemma~\ref{lem:facts}~(\ref{fact2}) above.

Conversely, suppose that there exists $\ell\in\{1,2,\ldots ,2k\}$ and
$\mathbf{c}\in H_\ell$ such that $\sum_{i=1}^mc_i\mathbf{v}_i=0$. By
Lemma~\ref{lem:facts}~(\ref{fact3}), we may write $\mathbf{c}$ as the
sum of $\ell$ difference vectors. Since multiplying a difference
vector by the scalar $-1$ produces another difference vector, we may
write $\mathbf{c}=\mathbf{a}-\mathbf{b}$, where
$\mathbf{a},\mathbf{b}$ are the sum of $\lfloor \ell/2\rfloor$ and
$\lceil \ell/2\rceil$ difference vectors respectively. Note that
$\mathbf{a}\neq\mathbf{b}$ since $\mathbf{c}\neq 0$. But
$\mathbf{a}\in H_{\lfloor \ell/2\rfloor}$ and $\mathbf{b}\in H_{\lceil
\ell/2\rceil}$, where $0\leq \lfloor \ell/2\rfloor\leq
\lceil\ell/2\rceil\leq \lceil 2k/2\rceil=k$, and so
Theorem~\ref{thm:khop_coverage} implies that $D$ does not meet the
bound, as required.
\end{IEEEproof}

\begin{definition}
Let $A$ be an abelian group. Let
$D=\{\mathbf{v}_1,\mathbf{v}_2,\ldots,\mathbf{v}_m\}\subseteq A$ be a
sequence of elements of $A$. We say that $D$ is a \emph{$B_h$ sequence
over $A$} if all the sums
\begin{equation}
\label{Bhdefn_eqn} \mathbf{v}_{i_1}+\mathbf{v}_{i_2}+\cdots +\mathbf{v}_{i_h}
\text{ with }1\leq i_1\leq\cdots\leq i_h\leq m
\end{equation}
are distinct.
\end{definition}

$B_h$ sequences (sometimes known as $B_h$-sets) have been studied
for many years, mainly in the case where $A=\mathbb{Z}$. See Graham~\cite{Grahamsurvey}, Halberstam
and Roth~\cite{HalberstamRoth66},
Lindstr\"{o}m~\cite{Lindstrom72}, O'Bryant~\cite{dynamicsurvey}, for example.
\begin{example}\label{bosechowlaconstruction}
Let $q$ be a prime power, let $h$ be an integer such that $h\geq 2$ and let $\alpha$ be a primitive element of
$\mathrm{GF}(q^h)$.  Bose and Chowla~\cite{BoseChowla62} have shown
that the set $\{a\in \mathbb{Z}_{q^h-1} \vert \alpha^a-\alpha\in
\mathrm{GF}(q)\}$ is a $B_h$ set in $\mathbb{Z}_{q^h-1}$ containing $q$
elements.
\end{example}
The following theorem demonstrates the relation between $B_h$ sequences and distinct difference configurations.
\begin{theorem}
\label{thm:DD_B} Let $k$ be a fixed integer, where $k\geq 2$.  Let
$D=\{\mathbf{v}_1,\mathbf{v}_2,\ldots,\mathbf{v}_m\}\subseteq
\mathbb{Z}^2$. Then $D$ is a $\DD(m)$ with maximal $k$-hop coverage if
and only if $D$ is a $B_{2k}$ sequence over $\bZ^2$.
\end{theorem}
\begin{IEEEproof}
Suppose $D$ is a $B_{2k}$ sequence over $\bZ^2$. We aim to show that $D$ is a
$\DD(m)$ with maximal $k$-hop coverage.

If $\mathbf{v}_i=\mathbf{v}_j$ for $i\not=j$ then $(2k-1)\mathbf{v}_1+\mathbf{v}_i=(2k-1)\mathbf{v}_1+\mathbf{v}_j$ and
so $D$ cannot be a $B_{2k}$ sequence. This contradiction implies
that the vectors are all distinct.

Suppose that $\mathbf{v}_i-\mathbf{v}_j=\mathbf{v}_{i'}-\mathbf{v}_{j'}$, where $i\not=j$,
$i'\not=j'$. Then $(2k-2)\mathbf{v}_1+\mathbf{v}_i+\mathbf{v}_{j'}=(2k-2)\mathbf{v}_1+\mathbf{v}_{i'}+\mathbf{v}_j$.
This contradicts the fact that $D$ is a $B_{2k}$ sequence, unless
$i=i'$ and $j'=j$. Thus $D$ has the distinct differences property.
Hence $D$ is a $\DD(m)$.

Suppose, for a contradiction, that $D$ does not have maximal $k$-hop
coverage. By Lemma~\ref{lem:maximalkhop_criterion} there exists
$\mathbf{a}=(a_1,a_2,\ldots,a_m) \in H_1\cup \cdots\cup H_{2k}$ such that
$\sum_{i=1}^ma_i\mathbf{v}_i=0$. Define $\mathbf{b}$ by
\[
b_i=\left\{\begin{array}{cl}
a_i&\text{when $a_i\geq 0$}\\
0&\text{ otherwise.}
\end{array}\right.
\]
Define $\mathbf{c}$ by the equation $\mathbf{a}=\mathbf{b}-\mathbf{c}$. Then the components of $\mathbf{b}$ and
$\mathbf{c}$ are all non-negative. Writing $t=\sum_{i=1}^m b_i=\sum_{i=1}^m
c_i=\sum_{a_i>0}a_i$, the definition of $H_1,H_2,\ldots,H_{2k}$
implies that $1\leq t\leq 2k$. Since $\mathbf{a}$ is non-zero, $\mathbf{b}\neq\mathbf{c}$.
But then our choice of $\mathbf{a}$ implies that
\[
(2k-t)\mathbf{v}_1+\sum_{i=1}^m b_i\mathbf{v}_i = (2k-t)\mathbf{v}_1 +\sum_{i=1}^m c_i\mathbf{v}_i.
\]
There are exactly $2k$ summands on both sides of this equality, so
$D$ cannot be a $B_{2k}$ sequence. This contradiction shows that
$D$ has maximal $k$-hop coverage, as required.

Now suppose that $D$ is a $DD(m)$ with maximal $k$-hop coverage.
Assume that $D$ is not a $B_{2k}$ sequence, so there exist two
distinct sums of the form~\eqref{Bhdefn_eqn} that are equal. By
cancelling terms that occur in both sums, we find that
$\sum_{i=1}^mb_i\mathbf{v}_i=\sum_{i=1}^mc_i\mathbf{v}_i$, where the
coefficients $b_i,c_i$ are all non-negative and where $\sum_{i=1}^m
b_i=\sum_{i=1}^m c_i=t$ for some integer $t$ such that $1\leq t\leq
2k$. But defining $a_i=b_i-c_i$ we find that $(a_1,a_2,\ldots,a_m)\in
H_t$ and $\sum_{i=1}^ma_i\mathbf{v}_i=0$. Hence $D$ does not have
maximal $k$-hop coverage, by Lemma~\ref{lem:maximalkhop_criterion}, as
required.
\end{IEEEproof}
The following construction converts a known construction for a
$B_{2k}$ sequence in $\mathbb{Z}_{q^{2k}-1}$ into a $B_{2k}$ sequence
in $\mathbb{Z}^2$, which is a $\DD(m)$ with maximal $k$-hop coverage
by Theorem~\ref{thm:DD_B}.
\begin{construction}
\label{BoseChowla_construction}
Let $k$ be a fixed integer such that $k\geq 2$. Let $q$ be a prime
power, and let $q^{2k}-1=ab$ where $a$ and $b$ are coprime. Then there
exists a set $X\subseteq \mathbb{Z}^2$ of dots that is doubly periodic
with periods $a$ and $b$, and such that the intersection of $X$ with
any $b\times a$ rectangle is a $\DD(q)$ with maximal $k$-hop coverage.
\end{construction}
\begin{IEEEproof}
The construction of Bose and Chowla~\cite{BoseChowla62} described in
Example~\ref{bosechowlaconstruction} shows there is a
$B_{2k}$ sequence over $\mathbb{Z}_{q^{2k}-1}$ consisting of $q$
elements. Note that by the Chinese Remainder Theorem there is a group
isomorphism $\mathbb{Z}_{q^{2k}-1}\rightarrow \mathbb{Z}_a\times
\mathbb{Z}_b$ given by $x\mapsto(x\bmod{a},x\bmod{b})$.  Thus there
are elements $\overline{\mathbf{v}}_1,\overline{\mathbf{v}}_2,\ldots
,\overline{\mathbf{v}}_q\in \mathbb{Z}_a\times \mathbb{Z}_b$ that form
a $B_{2k}$ sequence over $\mathbb{Z}_a\times \mathbb{Z}_b$. Let
$\rho:\mathbb{Z}^2\rightarrow \mathbb{Z}_a\times \mathbb{Z}_b$ be the
map defined by $\rho((x,y))=(x\bmod a,y\bmod b)$. We define
$X\subseteq \mathbb{Z}^2$ to be the set of vectors
$\mathbf{\mathbf{v}}\in \mathbb{Z}^2$ such that
$\rho(\mathbf{v})\in\{\overline{\mathbf{v}}_1,\overline{\mathbf{v}}_2,\ldots
,\overline{\mathbf{v}}_q\}$.

Since $\rho((x,y))=\rho((x+ia,y+jb))$ for any $i,j\in\mathbb{Z}$, we
see that $X$ is doubly periodic with periods $a$ and $b$
respectively. Let $R$ be an $b\times a$ rectangle in
$\mathbb{Z}^2$. For all $i\in\{1,2,\ldots ,m\}$, there is a unique
$\mathbf{v}_i\in R$ such that
$\rho(\mathbf{v}_i)=\overline{\mathbf{v}}_i$. Hence $X\cap
R=\{\mathbf{v}_1,\mathbf{v}_2,\ldots \mathbf{v}_q\}$. Moreover,
$\mathbf{v}_1,\mathbf{v}_2,\ldots,\mathbf{v}_q$ form a
$B_{2k}$ sequence over $\mathbb{Z}^2$, since if there are two sums of
the form~\eqref{Bhdefn_eqn} that are equal, then the images of these
sums under $\rho$ are also equal, which contradicts the fact that
$\overline{\mathbf{v}}_1,\overline{\mathbf{v}}_2,\ldots
,\overline{\mathbf{v}}_q$ form a $B_{2k}$ sequence over
$\bZ_a\times\bZ_b$. Thus
$\mathbf{v}_1,\mathbf{v}_2,\ldots,\mathbf{v}_q$ form a $\DD(q)$ with
maximal $k$-hop coverage by Theorem~\ref{thm:DD_B}, as required.
\end{IEEEproof}

This construction can be used to prove the existence of a $\DD(m,r)$
with maximal $k$-hop coverage where $r$ is small:

\begin{theorem}\label{thm:max_k_hop}
Let $k$ be a fixed integer such that $k\geq 2$. Define
$c=(\pi/16)2^{1/k}$.  Then there exists a $\DD(m,r)$ with maximal
$k$-hop coverage such that $m\sim cr^{1/k}$.
\end{theorem}

\begin{IEEEproof}
Let ${\cal S}\subseteq \mathbb{Z}^2$ be the set of points in
$\mathbb{Z}^2$ contained in a circle of radius $\lfloor r/2\rfloor$
about the origin. Note that $|{\cal S}|=(\pi/4) r^2+O(r)$ (by the
Gauss Circle Problem).

Let $q$ be the smallest prime power such that $q^{k}> 2r$. We have
that $q\leq (2r)^{1/k}+((2r)^{1/k})^{5/8}$ whenever $r$ is
sufficiently large by a classical result of Ingham~\cite{Ingham37} on
the gaps between primes. In particular, $q\sim (2r)^{1/k}$.

Define the integer $a$ by
\[
a=\left\{\begin{array}{cl} q^k-1&\mbox{ when $q$ is even,}\\
(q^k-1)/2&\mbox{ when $q^k\equiv 3\bmod 4$,}\\ (q^k+1)/2&\mbox{ when
$q^k\equiv 1\bmod 4$.}
\end{array}\right.
\]
Define $b=(q^{2k}-1)/a$. Since $\gcd(q^k-1,q^k+1)= 1$ when $q$ is even
and $\gcd(q^k-1,q^k+1)= 2$ when $q$ is odd, we find that $a$ and $b$
are coprime. Moreover, our choice of $q$ shows that $r\leq a\leq
b$. Let $X$ be the set of dots in $\mathbb{Z}^2$ given in
Construction~\ref{BoseChowla_construction}.

The average number of dots in a shift of $S$ by an element of
$\mathbb{Z}^2$ is $|S|q/(ab)$, and so we can find a shift $T$ of $S$
such that $|T\cap X|\geq |{\cal S}| q/(ab)$. Define $D\subseteq T\cap X$ to
be a subset of size $m$, where $m=\lceil |{\cal S}|q/(ab)\rceil$. Note that
$m \sim (\pi/4) r^2 q/(2r)^2\sim (\pi/16)2^{1/k}r^{1/k}$. Since $T$ is
a sphere of radius $\lfloor r/2\rfloor$, any pair of dots in $D$ are
at distance at most $r$. Moreover, the fact that $r\leq a\leq b$
implies that $T$ is contained in a $b\times a$ rectangle $R$. By
Construction~\ref{BoseChowla_construction}, $R\cap X$ is a $\DD(q)$
with maximal $k$-hop coverage. Since $D\subseteq T\cap X\subseteq
R\cap S$, we see that $D$ is a $\DD(m,r)$ with maximal $k$-hop
coverage. So the theorem follows, as required.
\end{IEEEproof}
Combining Theorems~\ref{thm:hextosquare},~\ref{thm:khopsame}
and~\ref{thm:max_k_hop}, we have the analogous result for the
hexagonal grid:
\begin{corollary}
\label{cor:hex_max_k_hop}
Let $k$ be a fixed integer such that $k\geq 2$. Define
$c^\prime=(\pi/16)2^{1/k}\left(\frac{2}{3}\right)^{1/2k}$. Then there
exists a $\DD^*(m,r)$ with maximal $k$-hop coverage such that $m\sim
c^\prime r^{1/k}$.
\end{corollary}

For any fixed values of $m$ and $k$, we define $r(k,m)$ to be the
smallest value of $r$ such that there exists a $\DD(m,r)$ with maximal
$k$-hop coverage. It is an important problem to determine
$r(k,m)$. The construction in Theorem~\ref{thm:max_k_hop} provides an
upper bound on $r(k,m)$, showing that when $k$ is fixed and
$m\rightarrow\infty$ we have $r(k,m)=O(m^k)$. We now provide a
corresponding lower bound on $r(k,m)$, which shows that the
construction in Theorem~\ref{thm:max_k_hop} is reasonable:
\begin{theorem}
\label{thm:bCmk} Let $k$ be an integer such that $k \geq 2$.
Then $\frac{m^k}{\sqrt{\pi}k!\cdot k}+o(m^k) \leq r(k,m) \leq
\frac{1}{2}\left(\frac{16}{\pi}\right)^km^k+ o(m^k)$.
\end{theorem}
\begin{IEEEproof}
The upper bound is proved in Theorem~\ref{thm:max_k_hop}.

To prove the lower bound, let $D$ be a $\DD(m,r)$ with maximal $k$-hop
coverage, where $r=r(k,m)$. The definition of maximal $k$-hop coverage
and Theorem~\ref{thm:khop_character} show that
$C_k(D)=\sum_{i=1}^k|H_i|$.  Let $B= \{ (a_1,a_2,\ldots,a_m) \in H_k
~:~ | \{ i : a_i \neq 0 \} |=2k \}$. Clearly $| B | =
\frac{m!}{(m-2k)!  k!^2}$ and
\[
\sum_{i=1}^k | H_i | =
\frac{m!}{(m-2k)! k!^2} +o(m^{2k})=\frac{m^{2k}}{k!^2} +o(m^{2k}).
\]
So $C_k(D)=\frac{m^{2k}}{k!^2} +o(m^{2k})$.

Every vector counted by $C_k(D)$ is the sum of at most $k$ difference
vectors of $D$. Each difference vector has length at most $r$, and so
every vector counted by $C_k(D)$ is contained in a circle of radius
$kr$ centred at the origin. Such a circle contains at most $\pi
(kr)^2+O(r)$ vectors in $\bZ^2$ (by Gauss's solution to the Gauss
circle problem). Thus
\[
\frac{m^{2k}}{k!^2} +o(m^{2k})= C_k(D)\leq \pi (kr)^2+O(r),
\]
which implies the lower bound of the theorem, as required.
\end{IEEEproof}
For the hexagonal grid, we denote the smallest $r$ for which there
exists a $\DD^*(m,r)$ with complete $k$-hop coverage by $r^*(m,k)$.
Combining Theorems~\ref{thm:bCmk} and \ref{thm:hextosquare}, we have
the following.
\begin{theorem}
\label{thm:rstar}
If $k \geq 2$ then $\sqrt{\frac{3}{2}}\frac{m^k}{\sqrt{\pi}k!\cdot k}+o(m^k)
\leq r^*(k,m) \leq \sqrt{\frac{3}{2}}\,\frac{1}{2}\left(\frac{16}{\pi}\right)^km^k+
o(m^k)$.
\end{theorem}
In the case $k=1$, we can use the results of \cite{BEMP1} to give
tighter bounds, as every distinct difference configuration has a
one-hop coverage of $m(m-1)$, which is thus maximal.
\begin{theorem}
\label{thm:bound_r1m} We have that
\[
\frac{2}{\sqrt{\pi}} m + o(m) \leq r(1,m) \leq
\frac{2}{\mu}m+o(m),
\]
where $\mu\approx 0.914769$ is the maximum value of
$((\pi/2)-2\theta+\sin 2\theta)/\cos\theta$ on the interval $0\leq
\theta\leq\pi/4$.
\end{theorem}
\begin{IEEEproof}
It is proved in~\cite{BEMP1} that if a $\DD(m,r)$ exists, then $m\leq
\frac{\sqrt{\pi}}{2}r+O(r^{2/3})$, which gives rise to the lower bound
on $r(1,m)$.  Furthermore, \cite{BEMP1} contains a construction of a
$\DD(m,r)$ with $m=(\mu/2)r+o(r)$ dots, from which we derive the upper
bound.
\end{IEEEproof}
The paper~\cite{BEMP1} also contains analogous results in the
hexagonal grid. From these, we can deduce the following bounds on
$r^*(1,m)$:
\begin{theorem} We have that
\[
\frac{\sqrt{2}\,3^{1/4}}{\sqrt{\pi}} m + o(m) \leq r^*(1,m) \leq
\frac{2^{1/2}3^{1/4}}{\mu}m+o(m),
\]
where $\mu$ is defined as in Theorem~\ref{thm:bound_r1m}.
\end{theorem}

Recall that we introduced the Manhattan and hexagonal metrics on the
square and hexagonal grids respectively in
Section~\ref{sec:models}. We conclude this subsection with a brief
discussion about the situation when we use these metrics rather than
Euclidean distance. For integers $k$ and $m$, define
$\overline{r}(k,m)$ to be the smallest integer $r$ such that there
exists a $\oDD(m,r)$ with maximal $k$-hop coverage, and define
$\overline{r}^*(k,m)$ to be the smallest integer $r$ such that there
exists a $\oDD^*(m,r)$ with maximal $k$-hop coverage.

\begin{theorem}
Let $k$ be a fixed integer, $k\geq 2$.  There exist constants
$c_1,c_2,c_3$ and $c_4$ such that for all sufficiently large integers $m$
\begin{align*}
c_1 m^k&\leq \overline{r}(k,m)\leq c_2m^k\text{ and}\\
c_3 m^k&\leq \overline{r}^*(k,m)\leq c_4m^k.
\end{align*}
\end{theorem}
\begin{IEEEproof}
By Theorem~\ref{euclee}, a $\oDD(m,r)$ with maximal $k$-hop coverage is
also a $\DD(m,r)$ with maximal $k$-hop coverage. So $r(k,m)\leq
\overline{r}(k,m)$. Moreover, a $\DD(m,r)$ with maximal $k$-hop
coverage is a $\oDD(m,\lceil \sqrt{2}r\rceil)$ with maximal $k$-hop
coverage, so $\overline{r}(k,m)\leq \lceil \sqrt{2}\,r(k,m)\rceil$. The
first statement of the theorem now follows by Theorem~\ref{thm:bCmk}.

The proof of the second statement of the theorem is similar, using
Theorems~\ref{thm:euchex} and~\ref{thm:rstar} in place of
Theorems~\ref{euclee} and~\ref{thm:bCmk} respectively.
\end{IEEEproof}

The results in~\cite{BEMP1} can be used to establish the following:

\begin{theorem}
We have that
\[
\overline{r}(1,m) = \sqrt{2}\,m+o(m).
\]
Moreover,
\[
(2/\sqrt{3})m+o(m) \leq\overline{r}^*(1,m)\leq (2/\mu)m+o(m),
\]
where $\mu=(2/3)^{3/2}(1+2\sqrt{7})/(\sqrt{2+\sqrt{7}})\approx 1.58887$.
\end{theorem}

\subsection{Minimum $k$-hop coverage}\label{sec:min}

Having established an upper bound for the $k$-hop coverage of a
$\DD(m)$ (and hence of a $\DD^*(m)$), we now consider the smallest
values it can take.
\begin{theorem}\label{easybound}
The $k$-hop coverage of a $\DD(m)$ is at least $km(m-1)$.
\end{theorem}
\begin{IEEEproof}
The one-hop coverage of a $\DD(m)$ is $m(m-1)$.

For $D=\{\mathbf{v}_1,\mathbf{v}_2,\dotsc,\mathbf{v}_m\}$ a $\DD(m)$,
let $\mathbf{u}=(d,e)$ be the difference vector with $\lvert d\rvert$
as large as possible. If there is more than one choice for
$\mathbf{u}$, choose $\mathbf{u}$ with $\lvert e\rvert$ as large as
possible subject to $\lvert d\rvert$ being maximal. Without loss of
generality, we can assume that $d >0$ and $e \geq 0$ (if not we can
flip and rotate the array to obtain an equivalent array with such
vector).

Let $S_1$ be the set of $m(m-1)$ vectors that can be reached by
one-hop paths from the origin. Then $S_1$ can be written as the
disjoint union of the two sets
\begin{equation*}
S_1^+=\{(x,y) \vert (x,y) \in S_1 ,~ x
>0 ~ \text{or} ~ (x=0 ~ \text{and} ~ y >0) \}\\
\end{equation*}
and $S_1^-=\{-(x,y) \vert(x,y) \in S_1^+ \}$.

For $i>1$, we define
\begin{align*}
S_i &= \{\mathbf{w} + (i-1)\mathbf{u} \vert \mathbf{w} \in S_1^+\}\,\,\cup\\
& \quad\quad\{(-\mathbf{w} - (i-1)\mathbf{u} \vert\mathbf{w} \in S_1^+\}.
\end{align*}
As $\bf{u}$ is a difference vector of $D$, the vectors of $S_i$ can
all be reached by $i$-hop paths from the origin.  Furthermore, $S_i
\cap S_j = \emptyset$ for $i \neq j$ and $|S_i|=m(m-1)$. Hence, the
theorem is proved.
\end{IEEEproof}

For certain values of $m$ there exist $\DD(m)$ for which the above
bound is tight.  For example, consider the following $\DD(3)$:
{\tiny \begin{equation*}{ \setlength{\arraycolsep}{.4\arraycolsep}
\begin{array}{c|c|c|c}

&&&\\\hline
\bullet&\bullet&\phantom{\bullet}&\bullet\\\hline
&&&\end{array}}
\end{equation*}}
The difference vectors in this example are
$\{\pm(1,0),\pm(2,0),\pm(3,0)\}$, and hence any of the $6k$ vectors of
the form $\pm(t,0)$ for $0<t\leq 3k$ can be reached by a $k$-hop path.

We can construct more examples where the bound is tight as follows. A
\emph{Golomb ruler} is a set $M$ of $m$ integers such that the differences
$x-y$ where $x,y\in M$ and $x\not=y$ are all distinct. A Golomb ruler is \emph{perfect} if
\begin{equation*}
\{u-v:u,v\in S\}=\{i\in\bZ: |i|\leq m(m-1)/2\}.
\end{equation*}
For example, the sequence $\{0,1,3\}$ is a perfect Golomb ruler.  The
$\DD(3)$ above was constructed from this sequence by taking
appropriate multiples of the vector $(1,0)$.  More generally, if $M$
is a perfect Golomb ruler then a configuration $D$ consisting of the
vectors $\mathbf{r}+i\mathbf{s}$ where $i\in M$ is a $\DD(m)$ with a
$k$-hop coverage of $km(m-1)$, and so meets the bound of
Theorem~\ref{easybound}. We say that $D$ is \emph{equivalent to a
perfect Golomb ruler} if we can construct it in this way. In fact, we
will now show that a $\DD(m)$ meets the bound of
Theorem~\ref{easybound} if and only if it is equivalent to a perfect
Golomb ruler.

\begin{lemma}
\label{lem:paralleldots}
Let $k$ be an integer, $k\geq 2$. Suppose $D$ is a $\DD(m)$ in which
there are differences $\mathbf{d}$ and $\mathbf{d^\prime}$ that are
not parallel. Then the $k$-hop coverage of $D$ is strictly greater
than $km(m-1)$.
\end{lemma}
\begin{IEEEproof}
Define the difference vector $\mathbf{u}$ and the sets $S_i$ as in the
proof of Theorem~\ref{easybound}. The set of difference vectors not
parallel to $\mathbf{u}$ is non-empty by assumption. Let $\mathbf{v}$
be a difference vector whose projection in the direction perpendicular
to $\mathbf{u}$ has length $p(\mathbf{v})$ as large as possible. Since
$k\geq 2$, the $k$-hop coverage of $D$ is at least
\[
\left|S_1\cup S_2\cup\cdots \cup S_k\cup\{2\mathbf{v}\}\right|.
\]
The argument in Theorem~\ref{easybound} shows the sets $S_i$ are
disjoint and have order $m(m-1)$. So the theorem follows if we can
show that $2\mathbf{v}\not\in S_1\cup S_2\cup\cdots \cup S_k$. But any vector in $S_i$ can be written in the form $\mathbf{w}\pm (i-1)\mathbf{u}$ where $\mathbf{w}$ is a difference vector, and therefore
\[
p(\mathbf{w}\pm (i-1)\mathbf{u})=p(\mathbf{w})\leq p(\mathbf{v})<2p(\mathbf{v})=p(2\mathbf{v}).
\]
Hence $2\mathbf{v}$ does not lie in any of the sets $S_i$, as required.
\end{IEEEproof}
\begin{theorem}
\label{thm:perfect}
Let $k$ be an integer such that $k\geq 2$, and let $D$ be a
$\DD(m)$. Then $D$ meets the bound of Theorem~\ref{easybound} if and
only if it is equivalent to a perfect Golomb ruler.
\end{theorem}
\begin{IEEEproof}
It is easy to see that if $D$ is equivalent to a perfect Golomb ruler,
then $D$ meets the bound of Theorem~\ref{easybound}.

Let $D$ be a $\DD(m)$ that meets the bound of
Theorem~\ref{easybound}. The set $S_\ell$ defined in the proof of
Theorem~\ref{easybound} is a set of $m(m-1)$ vectors that can be
reached by an $\ell$-hop path from the origin, but cannot be reached
by a path of length $\ell-1$. Thus $C_k(D)\geq C_{2}(D)+(k-2)m(m-1)$,
so $D$ meets the bound of Theorem~\ref{easybound} in the case
$k=2$. So to prove the theorem, we need only consider the case $k=2$.

Let $\mathbf{r}$ be a vector in
$D$. Lemma~\ref{lem:paralleldots} implies that all the difference
vectors in $D$ are parallel to a fixed vector $\mathbf{u}$. Let
$\mathbf{s}$ be the shortest vector in $\bZ^2$ that is parallel to
$\mathbf{u}$. Then (since $\bZ^2$ is a lattice) $D\subseteq
\{\mathbf{r}+i\mathbf{s}\mid i\in\bZ\}$. Thus $D$ is equivalent to a
Golomb ruler $M\subseteq\bZ$. Without loss of
generality, we may assume that the greatest common divisor of the
elements of $M$ is $1$, for if the greatest common divisor is $a$ then
we can replace $\mathbf{s}$ by $a\mathbf{s}$ and $M$ by $(1/a)M$.

It remains to show that $M$ is perfect. The set $S=\{x-y\vert x,y\in
M\}$ contains $m(m-1)+1$ elements, since $M$ is a Golomb ruler.  A
square reachable from the origin by a one-hop or two-hop path
corresponds to an element of $S+S=\{a+b\vert a,b\in S\}$. It is a
well-known result of additive combinatorics that for a set $A$ of
integers with $\lvert A \rvert=n$ it holds that $\lvert
A+A\rvert=2n-1$ if and only if the elements of $A$ are in arithmetic
progression.  The bound of Theorem~\ref{easybound} requires $S+S$ to
have size $2m(m-1)+1$ (due to the inclusion of 0); as this is equal to
$2\lvert S \rvert-1$ it follows that the elements of $S$ are in
arithmetic progression. Since $S=-S$ and the greatest common divisor
of the elements of $M$ is $1$ we find that $S=\{x\in\bZ\mid
|x|\leq m(m-1)/2\}$. So $M$ is a perfect Golomb ruler, as
required.
\end{IEEEproof}

\section{A $\DD(m)$ with Complete Two-Hop Coverage in a Rectangle}
\label{sec:complete-two-hop}

In Section~\ref{sec:twohop} we explored the range of values that the $k$-hop
coverage of a distinct difference configuration can take. When choosing a distinct difference configuration for use in Scheme~\ref{BEMPscheme} it may seem desirable to select a configuration with maximal two-hop coverage. However, from Theorem~\ref{thm:bCmk} we see that a $\DD(m,r)$ with maximal two-hop coverage has ``approximately'' $m^2=r$, which places too great a restriction on the maximum number of keys that each node can store in the resulting scheme. From a practical perspective it thus may be desirable to focus on connectivity within a localised region.

In this section we give a
construction of a $\DD(m)$ that ensures a two-hop path between a given
point $\mathbf{x}$ and any other grid point within a
$(2p-3)\times(2p-1)$ rectangle centred at $\mathbf{x}$, where $p$ is any prime greater than or equal to five. This allows the region to be tailored to the requirements of a specific application environment.

Our
construction can be thought of as being based on the periodicity
properties of a $B_2$ sequence in $\mathbb{Z}_{(p^2-p)}$ proposed by
Ruzsa in \cite{ruzsa}, or as a consequence of a periodic
generalisation of the Welch construction of a Costas array
\cite{GoTa84}.  In Subsection~\ref{sec:wpa} we discuss some
properties of a related doubly periodic array that we will exploit
later. In Subsection~\ref{sec:LocalDDM} we present the construction and demonstrate that it achieves complete two-hop coverage.

\subsection{The Welch Periodic Array}\label{sec:wpa}

\begin{definition} {\bf (Welch Periodic Array)} Let $\alpha$ be a
primitive root modulo a prime $p$.  We define the {\em
Welch periodic array} to be the set
\begin{equation*}
{\cal R}_p=\{(i,j)\in\mathbb{Z}^2\vert \alpha^j\equiv i \bmod p\}.
\end{equation*}
\end{definition}

This array is doubly periodic in the sense that if ${\cal R}_p$
contains a dot at position $(i,j)$ then it also contains dots at all
positions of the form $(i+\lambda p,j+\mu (p-1))$ where
$\lambda,\mu\in\mathbb{Z}$.  It has a distinct difference property
``up to periodicity'': see the lemma below. We say that
dots $A$ and $A'$ at positions $(i,j)$ and $(i',j')$ are
\emph{equivalent}, and we write $A\equiv A'$, if $i'=i+\lambda p$ and
$j'=j+\mu(p-1)$ for some $\lambda,\mu\in\bZ$.

\begin{lemma}\label{lem:ddperiod}
Let $d$ and $e$ be integers such that $d\not\equiv 0\bmod p$ and
$e\not\equiv 0\bmod (p-1)$. Suppose that ${\cal R}_p$ contains dots
$A$ and $B$ at positions $(i_1,j_1)$ and $(i_1+d,j_1+e)$ respectively,
and dots $A'$ and $B'$ at positions $(i_2,j_2)$ and $(i_2+d,j_2+e)$
respectively. Then $A\equiv A'$ and $B\equiv B'$.
\end{lemma}
\begin{IEEEproof}
By the definition of ${\cal R}_p$ we have
\begin{align*}
i_1&\equiv \alpha^{j_1} \bmod p\\
i_2&\equiv \alpha^{j_2} \bmod p\\
i_1+d&\equiv \alpha^{j_1+e} \bmod p\\
i_2+d&\equiv \alpha^{j_2+e} \bmod p.\\
\end{align*}
Eliminating $i_1$, $i_2$ and $d$ from these equations we get
\begin{equation*}
(\alpha^e-1)(\alpha^{j_1}-\alpha^{j_2})\equiv 0\bmod p.
\end{equation*}
Since $e\not\equiv 0\bmod (p-1)$, this implies that $j_1\equiv
j_2\bmod (p-1)$. The first two equations above then imply that
$i_1\equiv i_2\bmod p$.
\end{IEEEproof}
We note that in addition, if ${\cal R}_p$ contains dots at $(i,j)$ and
$(i+d,j)$ then $d\equiv 0 \bmod p$ and if it contains dots at $(i,j)$
and $(i,j+e)$ then $e\equiv 0 \bmod (p-1)$.  Thus we see that a vector
$(d,e)$ can occur at most once as a difference between two of the dots
of ${\cal R}_p$ that lie within any particular $(p-1)\times p$
rectangle.

\subsection{Construction of the $\DD(m)$}\label{sec:LocalDDM}

We now define a $\DD(m)$ by choosing a finite subset of the dots in
${\cal R}_p$, as follows.
\begin{construction}\label{con:complete}
Let $p$ be an odd prime. Let
$(i,j)\in\bZ^2$ be such that ${\cal R}_p$ has dots at $(i,j)$ and
$(i+1,j+1)$. Note that such a position $(i,j)$ exists. To see this, let
$i$ and $j$ be integers such that
\[
\alpha^j\equiv i\equiv \frac{1}{\alpha-1} \bmod p.
\]
The right-hand side of this equality is well-defined and non-zero
modulo $p$, and so there is a suitable choice for $i$ and $j$. Clearly
${\cal R}_p$ has a dot at the position $(i,j)$. But there is also a dot at
$(i+1,j+1)$ since
\[
\alpha^{j+1}\equiv \frac{\alpha}{\alpha-1}\equiv
\frac{1}{\alpha-1}+1\equiv i+1\bmod p.
\]
Consider the $(p-1)\times p$ rectangle $S$ bounded by the positions
$(i,j)$, $(i+p-1,j)$, $(i,j+p-2)$ and $(i+p-1,j+p-2)$.  By
construction, ${\cal R}_p$ has $p-1$ dots in $S$.  Due to its periodic
nature, ${\cal R}_p$ also has dots at positions $(i,j+(p-1))$,
$(i+p,j)$ and $(i+p+1,j+p)$.  We construct a configuration $\cal B$ by
adding these three dots to the set of dots in ${\cal R}_p\cap S$.
\end{construction}

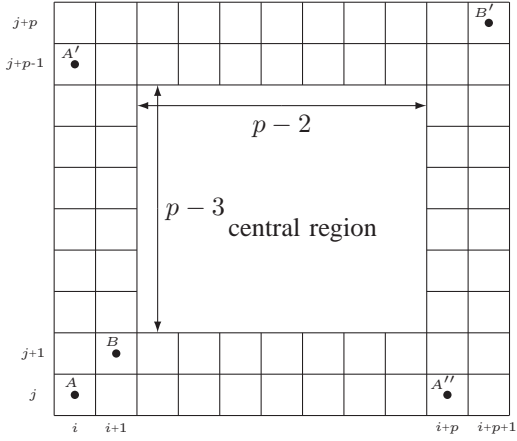
\begin{figure}[t]
\centering
\setlength{\unitlength}{1.1mm}
\begin{picture}(63,56)(-2,0)
\put(5,5){\line(1,0){55}}
\put(5,10){\line(1,0){55}}
\put(5,15){\line(1,0){55}}
\put(5,45){\line(1,0){55}}
\put(5,50){\line(1,0){55}}
\put(5,55){\line(1,0){55}}
\put(5,5){\line(0,1){50}}
\put(10,5){\line(0,1){50}}
\put(15,5){\line(0,1){50}}
\put(50,5){\line(0,1){50}}
\put(55,5){\line(0,1){50}}
\put(60,5){\line(0,1){50}}
\multiput(20,5)(5,0){6}{\line(0,1){10}}
\multiput(20,45)(5,0){6}{\line(0,1){10}}
\multiput(5,20)(0,5){5}{\line(1,0){10}}
\multiput(50,20)(0,5){5}{\line(1,0){10}}
\put(-1,50){\makebox(5,5){\tiny$j$+$p$}}
\put(-1,45){\makebox(5,5){\tiny$j$+$p$-$1$}}
\put(0,5){\makebox(5,5){\tiny$j$}}
\put(0,10){\makebox(5,5){\tiny$j$+$1$}}
\put(5,1){\makebox(5,5){\tiny$i$}}
\put(10,1){\makebox(5,5){\tiny$i$+$1$}}
\put(50,1){\makebox(5,5){\tiny$i$+$p$}}
\put(56,1){\makebox(5,5){\tiny$i$+$p$+$1$}}
\put(17.5,30){\vector(0,1){15}}
\put(17.5,30){\vector(0,-1){15}}
\put(19,27.5){\makebox(5,5){ $p-3$}}
\put(32.5,42.5){\vector(1,0){17.5}}
\put(32.5,42.5){\vector(-1,0){17.5}}
\put(30,37.5){\makebox(5,5){$p-2$}}

\put(25,25){\makebox(20,5){central region}}

\put(4.5,8){\makebox(5,2){\tiny$A$}}
\put(7.5,7.5){\circle*{1}}
\put(9.5,13){\makebox(5,2){\tiny$B$}}
\put(12.5,12.5){\circle*{1}}
\put(4.5,48){\makebox(5,2){\tiny$A^\prime$}}
\put(7.5,47.5){\circle*{1}}
\put(49.5,8){\makebox(5,2){\tiny$A^{\prime\prime}$}}
\put(52.5,7.5){\circle*{1}}
\put(54.5,53){\makebox(5,2){\tiny$B^\prime$}}
\put(57.5,52.5){\circle*{1}}
\end{picture}
\caption{The configuration $\cal B$. The five dots shown are the dots
that lie the border of width $2$ of the $(p+1)\times(p+2)$ rectangle
containing the configuration.}\label{fig:calB}
\end{figure}

Our configuration $\cal B$ is shown in Fig.~\ref{fig:calB}. The
configuration is contained in a $(p+1)\times (p+2)$ rectangle. The
\emph{border region} of width $2$ contains exactly $5$ dots: $A,A',A'',B$
and $B'$. The \emph{central region} is a $(p-3)\times (p-2)$
rectangle. This region contains $p-3$ dots: one column is empty, but
every other column and every row contains exactly one dot. Note that
$A\equiv A'\equiv A''$ and $B\equiv B'$, but there are no other
equivalent pairs of dots in $\cal B$.

\begin{lemma}
The configuration $\cal B$ is a $\DD(p+2)$, all of whose points lie in a $(p+1)\times (p+2)$ rectangle.
\end{lemma}
\begin{IEEEproof}
We have already remarked that $\cal B$ contains $p+2$ dots, all lying
in a $(p+1)\times (p+2)$ rectangle. So it remains to show that $\cal
B$ satisfies the distinct differences property.

Suppose, for a contradiction, that $X$ and $Y$, and $X'$ and $Y'$, are
distinct pairs of dots in $\cal B$ with the same difference vector
$(d,e)$.

Suppose that $d\in\{0,-p,p\}$ or $e\in\{0,-(p-1),(p-1)\}$. A
difference vector between a dot in the central region of our
configuration and any other dot has $x-$ and $y-$coordinates of
absolute value at most $p-1$ or $p-2$ respectively. Moreover, a
central dot is the only dot in its row and column. So our assumption
implies that none of $X,X',Y,Y'$ can lie in the central region of our
configuration. But the $5\times 4$ ordered pairs of dots in the
border region all have distinct difference vectors, and so we have a
contradiction in this case.

So we may assume that $d\not\in\{0,-p,p\}$ and
$e\not\in\{0,-(p-1),(p-1)\}$. In particular, since all dots lie in a
$(p+1)\times (p+2)$ rectangle, we see that $d\not\equiv 0\bmod p$ and
$e\not\equiv 0\bmod (p-1)$. Lemma~\ref{lem:ddperiod} now implies that
$X\equiv X'$ and $Y\equiv Y$. If $X=X'$ then $Y=Y'$ which contradicts
the fact that our pairs of dots are distinct. Hence $X\neq X'$. The
fact that $X\equiv X'$ now implies that $X$ and $X'$ must lie in the
border of our configuration. A similar argument implies the same is
true for $Y$ and $Y'$. As in the paragraph above, we now have a
contradiction. Thus the lemma follows.
\end{IEEEproof}

Our aim is to show (Theorem~\ref{thm:complete_coverage}) that $\cal B$
achieves complete two-hop coverage on a $(2p-3)\times(2p-1)$ rectangle
relative to the central point of the rectangle.  In order to
demonstrate this, it is necessary to show that every vector $(d,e)$
with $\lvert d \rvert \leq p-1$ and $\lvert e \rvert\leq p-2$ can be
expressed as a two-hop path of difference vectors from $\cal B$.  The
following lemma proves this for the majority of such vectors
$(d,e)$.

\begin{lemma}\label{lem:allnonzero}
Any vector of the form $(d,e)$, where $d$ and $e$ are non-zero integers satisfying $\lvert d \rvert\leq p-1$ and $\lvert e\rvert\leq p-2$, can be expressed as the sum of two difference vectors from $\cal B$.
\end{lemma}
\begin{IEEEproof}
Consider the $(p-1)\times p$ rectangle $S$ defined in
Construction~\ref{con:complete}, and let $\cal A$ be the restriction
of ${\cal R}_p$ to the $(2p-2)\times 2p$ subarray whose lower leftmost
corner coincides with that of $S$.

We partition $\cal A$ into four $(p-1)\times p$
subarrays as follows:
\begin{equation*}
\left(\begin{array}{c|c}
{\cal D}_3 & {\cal D}_4\\ \hline
{\cal D}_1 & {\cal D}_2
\end{array}
\right)
\end{equation*}
The periodicity of ${\cal R}_p$ means that the set of dots of ${\cal
R}_p$ contained in each subarray is a translation of the set of dots
of ${\cal R}_p$ contained in ${\cal D}_1$. Moreover, since ${\cal
D}_1=S$, all the dots in ${\cal D}_1$ are contained in $\cal B$.

We claim that each of the vectors $(d,e)$ appears as the difference of
two points in $\cal A$. Since the negative of a difference vector is
always a difference vector, we may assume without loss of generality
that $d> 0$. Suppose that $e> 0$. There is a unique position
$(i',j')\in {\cal D}_1$ such that
\[
\alpha^{j'}\equiv i'\equiv \frac{d}{\alpha^e-1}\bmod p.
\]
It is easy to check, just as in Construction~\ref{con:complete}, that
${\cal R}_p$ has dots at $(i',j')$ and $(i'+d,j'+e)$. Since $d$ and
$e$ are both positive, $(i'+d,j'+e)$ lies in $\cal A$, and so our
claim follows in this case. The argument for the case when $e<0$ is
exactly the same, except now we choose $(i',j')\in {\cal D}_3$. So the
claim follows.

To prove the lemma, we need to show that each difference vector
$(d,e)$ can be written as the sum of two difference vectors of
$\cal B$.  This follows from the paragraph above and the following
observations:
\begin{itemize}
\item Any vector connecting two dots of ${\cal D}_1$ is a difference
vector of $\cal B$ by construction.
\item Due to the periodicity of ${\cal R}_p$, a vector connecting a
dot in ${\cal D}_1$ with a dot in ${\cal D}_3$ (or, similarly, a
dot in ${\cal D}_2$ with a dot in ${\cal D}_4$) can be expressed
as the sum of the vector $(0,p-1)$ (which occurs as a difference
between the dots $A$ and $A'$ in $\cal B$) and some other difference
vector of $\cal B$.
\item A vector connecting a dot in ${\cal D}_1$ with a
dot in ${\cal D}_2$ (or, similarly, a dot in ${\cal D}_3$ with a
dot in ${\cal D}_4$) can be expressed as the sum of the difference
vector $(p,0)$ (which occurs between $A$ and $A''$) and some other
difference vector of $\cal B$.
\item A vector connecting a dot in ${\cal D}_1$ with a dot in
${\cal D}_4$ is the sum of the difference vector $(p,p-1)$ (which
occurs between $B$ and $B'$) and some other difference vector of $\cal
B$.
\item A vector connecting a dot in ${\cal D}_3$ with a dot in
${\cal D}_2$ is the sum of the difference vector $(p,-(p-1))$ (which
occurs between $A'$ and $A''$) and some other difference vector of
$\cal B$.
\end{itemize}
\end{IEEEproof}

It remains to consider vectors that have a zero co-ordinate. We will
use the following lemma in our proof that such vectors all occur as
the sum of two difference vectors from $\cal B$.

\begin{lemma}
\label{lem:diff+add}
Let $t$ be a positive integer with $t\geq 3$.
Let $\cF$ be a set of integers satisfying the following properties:
\begin{enumerate}[(a)]
\item \label{cFhast}$\lvert{\cF}\rvert=t+1$,
\item ${\cF}\subset\{-(t-1),-(t-2),\dotsc,-1\}\cup\{1,2,\dotsc,t-1\}\cup\{t+1\}$,
\item $\{1,-(t-1),t+1\}\subset \cF$,
\item $\exists i \in \cF\setminus \{1,-(t-1),t+1\}$ with $i<0$,
\item \label{justone}if $i>0$ and $i\in \cF\setminus \{1,-(t-1),t+1\}$ then $i-t\notin \cF$.
\end{enumerate}
Then each positive integer $\gamma$ with $1\leq \gamma \leq t-1$ has a
representation of the form $\gamma=j-i$ where $i,j\in \cF$.
\end{lemma}
\begin{IEEEproof}
Since $\cF\setminus \{1,-(t-1),t+1\}$ contains $t-2$ elements,
(\ref{justone}) implies that $\cF$ must contain precisely one element
of each pair $\{i,i-t\}$ for $i=2,3,\dotsc,t-1$.  Suppose, for a
contradiction, that there exists a positive integer $\gamma\leq
t-1$ that cannot be expressed as the difference between two elements
of $\cF$.

Suppose that $\gamma >1$.  Since $1,t+1\in\cF$, our assumption implies
that $1-\gamma\notin \cF$ and $t+1-\gamma \notin \cF$.  But
$1-\gamma=(t+1-\gamma)-t$, hence one of these numbers must be
contained in $\cF$, which gives a contradiction in this case.

Suppose that $\gamma=1$. The assumption implies that $\cF$ does not
contain a pair of integers that differ by 1.  If $t$ is odd this
implies that $\cF\setminus\{t+1\}$ contains at most $(t-1)/2$ positive
integers, and at most $(t-1)/2$ negative integers, hence $\cF$
contains at most $(t-1)+1=t$ integers, which contradicts
(\ref{cFhast}).  If $t$ is even, then in order for the size of $\cF$
to be $t+1$, $\cF\setminus\{t+1\}$ must contain $t/2$ positive
integers, all of which are odd, and $t/2$ negative integers that are
also all odd.  This implies that for each positive odd integer $1<i<t$
we have that $i\in \cF$ and $i-t\in\cF$, which contradicts
(\ref{justone}). So the lemma follows.
\end{IEEEproof}
We can now combine these two lemmas to obtain our desired result:
\begin{theorem}
\label{thm:complete_coverage}
Let $p$ be a prime, $p\geq 5$. The distinct difference configuration
$\cal B$ achieves complete two-hop coverage on a $(2p-3)\times(2p-1)$
rectangle relative to the central point of the rectangle.
\end{theorem}
\begin{IEEEproof}
By Lemma~\ref{lem:allnonzero}, any vector $(d,e)$ from the centre of a
$(2p-3)\times(2p-1)$ rectangle to another point of the rectangle can
be expressed as the sum of two difference vectors of $\cal B$ if $d$
and $e$ are non-zero.

We now consider vectors of the form $(0,e)$ with $0<e\leq p-2$.  Such
a vector can be expressed as the sum of two difference vectors of
$\cal B$ if $\cal B$ has difference vectors of the form $(1,y')$ and
$(1,y)$ with $y'-y=e$.  The second coordinates of the set
of difference vectors of $\cal B$ of the form $(1,y)$ with $y\neq 0$
satisfy the conditions of Lemma~\ref{lem:diff+add} for $t=p-1$, since:
\begin{enumerate}[(a)]
\item The left-most column of the array contains two dots; all
other columns contain a single dot apart from a single central column
which is empty. So $\cal B$ has $p$ difference vectors of the form
$(1,y)$ with $y\neq 0$.
\item Except for the vector $(1,p)$, all difference vectors of $\cal B$ of the form $(1,y)$ with $y\neq 0$ satisfy $|y|\leq p-2$.
\item The vectors $(1,1)$, $(1,-(p-2))$ and $(1,p+1)$ are all
difference vectors of $\cal B$ (as they occur as differences between
dots in the border region of $\cal B$, see Fig.~\ref{fig:calB}).
\item The difference vectors of $\cal B$ of the form $(1,y)$ cannot
all satisfy $y>0$. This is obvious if the right-most central column
contains a dot. If this column is empty and $y$ is always
positive, then the remaining $(p-3)\times (p-3)$ central region must
contain dots along a lower-left to top-right diagonal. Since $p\geq
5$, two central dots have the difference vector $(1,1)$. Since dots
$A$ and $B$ also have this difference vector, the distinct difference
property is violated and so we have a contradiction, as required.
\item If $(1,y)$ with $y\not=1,p$ is a difference vector of $\cal B$
then $(1,y-(p-1))$ is not. For Lemma~\ref{lem:ddperiod} implies that
the dots involved must be equivalent, and so must be in the border
region of our construction.
\end{enumerate}
Lemma~\ref{lem:diff+add} now implies that any vector $(0,e)$ with
$0<e\leq p-2$ has an expression in the form $(0,e)=(1,y')+(-1,-y)$
where $(1,y')$ and $(1,y)$ are difference vectors of $\cal B$.  Vectors
of the form $(0,e)$ with $-(p-2)<e<0$ can be written as
$(1,y)+(-1,-y')$.

In a similar manner, we can show that the first coordinates of the
difference vectors of $\cal B$ of the form $(x,1)$ satisfy the
conditions of Lemma~\ref{lem:diff+add} with $t=p$, and hence any
vector of the form $(d,0)$ with $0<\lvert d\rvert\leq p-1$ can be
written as the sum of two difference vectors of $\cal B$.  Thus the
result is proved.
\end{IEEEproof}

We can thus apply the $DD(m)$ specified in Construction~\ref{con:complete} to Scheme~\ref{BEMPscheme} in order to establish a key predistribution scheme which guarantees two-hop paths between a node and all of its neighbours within a surrounding rectangular region. This provides a powerful notion of local connectivity in order to facilitate connectivity across the wider network. The resulting scheme is also highly configurable, since the value of $p$ can be adjusted in order to tradeoff storage against the size of the fully connected local region.

\section{Conclusion and Open Problems}
\label{sec:conclude}

In this paper we have studied properties of distinct difference
configurations, which can be used to design efficient key
predistribution schemes for wireless sensor networks based on grids.

In Section~\ref{sec:twohop} we explored the $k$-hop coverage of a $\DD(m,r)$. 
We characterised maximal $k$-hop coverage
in terms of $B_{2k}$ sequences over $\bZ^2$, and we used a known
construction of $B_{2k}$ sequences over $\bZ$ to produce a $\DD(m,r)$
with maximal $k$-hop coverage and of the order of $r^{1/k}$ dots. We
provided an argument that shows that the order of magnitude of the number
of dots is correct (by bounding the functions $r(k,m)$). These results indicate the range of achievable parameters, which in turn determine the connectivity properties of the resulting key predistribution schemes. It would be
interesting to find better bounds on the leading coefficient of
$r(k,m)$, and it would be worthwhile determining $r(k,m)$ precisely
for small values of $k$ and $m$. Similar comments hold for the
function $r^*(k,m)$, and for the analogous situations using the
Manhattan or hexagonal metric.

In Section~\ref{sec:complete-two-hop} we constructed a $\DD(m,r)$ with
complete $2$-hop coverage within a large rectangular region centred on
the origin. This $\DD(m,r)$ can be used to design key predistribution schemes with excellent local connectivity properties.
The area of the fully connected region is of the order of $m^2$. It would be interesting to investigate whether
there are any constructions that achieve complete two-hop coverage in significantly
larger rectangles. Constructions that are optimised with respect to two-hop coverage for
other natural regions, for example a circle of large radius, would also be of practical interest. A further open problem is whether there exist any good constructions, for any natural region, achieving complete $k$-hop coverage for $k\geq 3$.

%
%

\begin{IEEEbiographynophoto}{Simon R. Blackburn}
received his BSc in Mathematics from the University of Bristol in
1989, and his DPhil in Mathematics from the University of Oxford in
1992. Since then he has worked at Royal Holloway, University of London
as a Research Assistant (1992-95), an Advanced Fellow (1995-2000), a
Reader in Mathematics (2000-2003) and a Professor in Pure Mathematics
(2004-). He was Head of the Mathematics Department from 2004 to
2007. His research interests include cryptography, group theory, and
combinatorics with applications to computer science.
\end{IEEEbiographynophoto}

\begin{IEEEbiographynophoto}{Tuvi Etzion}
(M'89-SM'99-F'04) was born in Tel Aviv, Israel, in 1956. He received
the B.A., M.Sc., and D.Sc. degrees from the Technion - Israel
Institute of Technology, Haifa, Israel, in 1980, 1982, and 1984,
respectively.

From 1984 he held a position in the department of Computer Science at
the Technion, where he has a Professor position. During the years
1986-1987 he was Visiting Research Professor with the Department of
Electrical Engineering - Systems at the University of Southern
California, Los Angeles. During the summers of 1990 and 1991 he was
visiting Bellcore in Morristown, New Jersey.  During the years
1994-1996 he was a Visiting Research Fellow in the Computer Science
Department at Royal Holloway, University of London. He also had
several visits to the Coordinated Science Laboratory at University of
Illinois in Urbana-Champaign during the years 1995-1998, two visits to
HP Bristol during the summers of 1996, 2000, several visits to the
department of Electrical Engineering, University of California at San
Diego during the years 2000-2009, and to the Mathematics department at
Royal Holloway, University of London during the years 2007-2009. His
research interests include applications of discrete mathematics to
problems in computer science and information theory, coding theory,
and combinatorial designs.

Dr Etzion was an Associate Editor for Coding Theory for
the IEEE Transactions on Information Theory from 2006 till 2009.
\end{IEEEbiographynophoto}

\begin{IEEEbiographynophoto}{Keith M. Martin} joined the Information
Security Group at Royal Holloway, University of London as a lecturer
in January 2000. He received his BSc (Hons) in Mathematics from the
University of Glasgow in 1988 and a PhD from Royal Holloway in
1991. Between 1992 and 1996 he held a Research Fellowship in the
Department of Pure Mathematics at the University of Adelaide,
investigating mathematical modeling of cryptographic key distribution
problems. In 1996 he joined the COSIC research group of the Katholieke
Universiteit Leuven in Belgium where he was primarily involved in an
EU ACTS project concerning security for third generation mobile
communications. He has also held visiting positions at the University
of Wollongong, University of Adelaide and Macquarie
University. Keith's current research interests include cryptography,
key management and wireless sensor network security.

Prof.\ Martin is an Associate Editor for Complexity and Cryptography
for IEEE Transactions on Information Theory.
\end{IEEEbiographynophoto}

\begin{IEEEbiographynophoto}{Maura B. Paterson}
received a BSc from the University of Adelaide in 2002 and a PhD from
Royal Holloway, University of London in 2005.  She has worked as a
research assistant in the Information Security Group at Royal
Holloway, and is currently at the Department of Economics, Mathematics
and Statistics at Birkbeck, University of London.  Her research
interests include applications of combinatorics in information
security.
\end{IEEEbiographynophoto}


\begin{thebibliography}{10}

\bibitem{BEMP1} S.R. Blackburn, T. Etzion, K.M. Martin, and
M.B. Paterson, ``Two-dimensional patterns with distinct differences:
Constructions, bounds and maximal anticodes,'' \emph{preprint},
2008.

\bibitem{BEMP} S.R. Blackburn, T. Etzion, K.M. Martin, and
M.B. Paterson, ``Efficient key predistribution for grid-based
wireless sensor networks,'' in (R. Safavi-Naini, Ed) \emph{Proc. ICITS
2008}, Lecture Notes in Computer Science~5155, Springer-Verlag,
Berlin, pp. 54--69, 2008.


\bibitem{romer-design}
K.~{R\"o}mer and F.~Mattern.
\newblock The design space of wireless sensor networks.
\newblock {\em {IEEE} Wireless Communications Magazine}, 11(6):54--61, 2004.

\bibitem{cysurvey} S.A.~\c{C}amtepe and B. Yener, ``Key Distribution
Mechanisms for Wireless Sensor Networks: a Survey'', \emph{Rensselaer
Polytechnic Institute Tech.~Report} TR-05-07 March 2005.

\bibitem{keithframework} K.M.~Martin and M.B.~Paterson, ``An
Application-Oriented Framework for Wireless Sensor Network Key
Establishment'', \emph{Electron. Notes Theor. Comput. Sci.},
vol.\,192, pp.\,31--41, 2008.

\bibitem{newsurvey} Y.~Xiao, V.K.~Rayi, B.~Sun, X.~Du, F.~Hu and M.~Galloway,
``A survey of key management schemes in wireless sensor
networks'', \emph{Comput. Commun.}, vol.\,30, pp.\,2314--2341,
2007.

\bibitem{eschenauer02}
L.~Eschenauer and V.D.~Gligor ``A key-management scheme for distributed sensor networks'', \emph{CCS '02: Proc. of the 9th ACM Conference on Computer and Communications Security}, pp.\,41--47, 2002.

\bibitem{soilnet}
Institut f\"{u}r {C}hemie und {D}ynamik der {G}eosph\"{a}re
({I}{C}{G}),
  {F}orschungszentrum {J}\"{u}lich: {S}oil{N}et -- a {Z}igbee based soil
  moisture sensor network.
\newblock http://www.fz-juelich.de/icg/icg-4/index.php?index=739, 2008.


\bibitem{nectarine}
Integrated smart sensing systems.
\newblock http://dpi.projectforum.com/isss/11, 2007.

\bibitem{irrigate}
J.~McCulloch, P.~McCarthy, S.~M. Guru, W.~Peng, D.~Hugo, and
A.~Terhorst.
\newblock Wireless sensor network deployment for water use efficiency in
  irrigation.
\newblock In {\em REALWSN '08: Proceedings of the Workshop on Real-world
  Wireless Sensor Networks}, pages 46--50, New York, NY, USA, 2008. ACM.

\bibitem{CoSl93}
J.H. Conway and N.J.A. Sloane, \emph{Sphere Packings, Lattices,
and Groups,} New York: Springer-Verlag, 1993.

\bibitem{leestin07}
J.~Lee and D.R.~Stinson ``On the construction of practical key predistribution schemes for distributed sensor networks using combinatorial designs'', \emph{ACM Trans. Inf. Syst. Secur.}, vol.\,11(2), pp.\,1--35, 2008.

\bibitem{Du03apairwise}
W.~Du, J.~Ding, Y.S.~Han, P.K.~Varshney, J.Katz and A.Khalili ``A pairwise key pre-distribution scheme for wireless sensor networks'', \emph{ACM Trans. Inf. Syst. Secur.}, vol.\,8, pp.\,228--258, 2005.

\bibitem{Chan03randomkey}
H.~Chan, A.~Perrig and D.~Song ``Random key predistirbution schemes for sensor networks'' \emph{IEEE Sumposium on Security and Privacy}, pp.\,197--213, 2003.

\bibitem{expander}
S.A.~\c{C}amtepe, B. Yener and M. Yung ``Expander graph based key distribution mechanisms in wireless sensor networks'', \emph{ICC '06, IEEE International Conference on Communications}, pp.\,2262--2267, 2006.

\bibitem{liuningli05}
D.~Liu, P.~Ning and R.~Li "Establishing pairwise keys in distributed sensor networks'', \emph{ACM Trans. Inf. Syst. Secur.}, vol.\,8(1), pp.\,41--77, 2005.

\bibitem{Grahamsurvey} S.W. Graham, ``$B_h$ sequences'',
\emph{Analytic Number Theory, Vol.~1 (Allerton Park, IL, 1995)},
Birkhauser, Boston, pp.\,431--449, 1996.

\bibitem{HalberstamRoth66} H. Halberstam and K.F. Roth,
\emph{Sequences, Volume I}, London: OUP, 1966.

\bibitem{Lindstrom72} B. Lindstr\"om, ``On $\mathrm{B}_2$ sequences of
vectors'', \emph{J Combinatorial Theory}, vol.\,4, pp. 261--265, 1972.

\bibitem{dynamicsurvey} K. O'Bryant, ``A complete annotated bibliography of work related to {S}idon sequences'', \emph{Electron. J. Combin} Dynamic Survey 11, 2004.

\bibitem{BoseChowla62} R.C. Bose and S. Chowla, ``Theorems in the
additive theory of numbers'', \emph{Comment. Math. Helvet.} vol.\,37, pp. 141--147, 1962--63.

\bibitem{Ingham37} A.E. Ingham, ``On the difference between consecutive
  primes'', \emph{Quart. J. Math. Oxford (Old Series)}, vol.\,8, pp. 255--266, 1937.

\bibitem{ruzsa}
I.Z. Ruzsa ``Solving a linear equation in a set of integers'', \emph{Acta\
Arith.}, vol.\,65, pp.\,259--282, 1993.

\bibitem{GoTa84}
S.W. Golomb and H. Taylor ``Constructions and properties of
Costas arrays'', \emph{Proceedings of the IEEE}, vol.\,72,
pp.\,1143--1163, 1984.

\end{thebibliography}
\end{document}